\documentclass{amsproc}
\usepackage{amssymb,amsthm,amscd}
 
\newtheorem{theorem}{Theorem}[section]
\newtheorem{proposition}[theorem]{Proposition}
\newtheorem{lemma}[theorem]{Lemma}
\newtheorem{corollary}[theorem]{Corollary}

\theoremstyle{remark}
\newtheorem{remark}[theorem]{Remark}

\theoremstyle{definition}
\newtheorem*{definition}{Definition}
\numberwithin{equation}{section}

\def\vomega{{\vec{\omega}}}
\def\complex{\mathbb{C}}
\def\real{\mathbb{R}}
\def\integer{\mathbb{Z}}
\def\BBf{\mathbf{B}}
\def\supp{\mathrm{supp}}
\def\sphere{\mathbf{S}^{d-1}}
\def\cone{\mathbf{C}} 
\def\Fourier{\mathbb{F}}

\def\BB{\mathcal{C}}
\def\CC{\mathcal{C}}
\def\DD{\mathcal{D}}
\def\FF{\mathcal{F}}

\def\II{\mathcal{I}}

\def\KK{\mathcal{K}}
\def\LL{\mathcal{L}}
\def\MM{\mathcal{M}}

\def\PP{\mathcal{P}}
\def\QQ{\mathcal{Q}}

\def\SS{\mathcal{S}}
\def\TT{\mathcal{T}}
\def\UU{\mathcal{U}}
\def\TTT{\mathbb{T}}
\def\VV{\mathcal{V}}
\def\WW{\mathcal{W}}
\def\tr{\mathrm{tr}}
\def\id{\mathrm{Id}}

\def\cc{{\Subset}}
\def\cB{\mathcal{B}}
\def\hB{\hat{\mathcal{B}}}
\def\hBBf{{\widehat{\BBf}}}

\begin{document}

\title[Dynamical determinants and spectrum] {Dynamical determinants 
and spectrum\\ for hyperbolic diffeomorphisms}

\author{Viviane Baladi}

\address{CNRS-UMR 7586, Institut de Math\'e\-ma\-ti\-ques  Jussieu, Paris, France}
\curraddr{D.M.A., UMR 8553,\'Ecole Normale Sup\'erieure,  75005 Paris, France}
\email{viviane.baladi@ens.fr}
\thanks{Part of this work was done during the trimester ``Time at Work"
in Institut Henri Poincar\'e, Paris, 2005.  VB also acknowledges useful comments from A. Avila and S. Gou\" ezel, support
from the French-Brazilian agreement, and hospitality of IMPA, Rio de Janeiro
and the Fields Institute, Toronto. MT is in part supported by Grant-in-Aid for Scientific Research (B) 18340044 from Japan Society for the Promotion of Science.}

\author{Masato Tsujii}
\address{Department of Mathematics, Kyushu University, 
Fukuoka,  Japan}

\email{tsujii@math.kyushu-u.ac.jp}

\dedicatory{Dedicated to Prof. Michael Brin on the occasion of his 60th birthday.}

\subjclass[2000]{37C30, 37D20,37A25}

\date{\today}

\begin{abstract}
For smooth  hyperbolic dynamical systems  and smooth weights,
we relate  Ruelle transfer operators with
dynamical Fredholm determinants and dynamical zeta
functions:
First, we establish bounds for the essential spectral radii 
of the transfer operator
on new spaces of anisotropic distributions,
improving our previous results \cite{BT}.
Then we give a new proof of Kitaev's \cite{Ki} lower bound for the radius of
convergence of the dynamical Fredholm determinant. 
In addition we show that  the zeroes of the determinant in the corresponding disc 
are in bijection with
the  eigenvalues of the transfer operator on our spaces of anisotropic distributions,
closing a question which remained open for a decade.
\end{abstract}

\maketitle

\section{Introduction}\label{sec:intro}

\subsection{Historical perspective}

The spectral properties of transfer operators and their relations to analytic properties of dynamical Fredholm determinants and dynamical zeta functions are fascinating subjects in study of smooth dynamical systems. The basic idea about the relation is rather simple: 
The dynamical Fredholm determinant of a transfer operator $\mathcal{L}$ 
associated to a dynamical system $T$ and a weight
function $g$ is formally defined by
\[
d_\LL(z)=\det({\mathrm Id}-z \mathcal{L}):=\exp \left(- \sum_{m=1}^\infty\frac {z^m}{m}
\sum_{T^m(x)=x}
\frac{g^{(m)}(x)} {|\det (\id-DT^m(x))|} \right) \, .
\]
Naturally, we would like that the inverse of each eigenvalue of the transfer operator
$\mathcal{L}$ appears as a zero of the dynamical Fredholm determinant. To make mathematically rigorous
statements, we first have to show that the transfer operator has nice spectral properties (similar
to those of compact operators) on a suitable Banach space. 
Finding the right Banach space is thus one of the hurdles in this subject. Then, we have to give an interpretation of the
sums over periodic points as (approximate) 
traces of transfer operators, recalling the formal relation
$$
\det({\mathrm{Id}}-z F)=\exp\left(-\sum_{m=0}^\infty 
\frac{z^m}{m} \mathrm{Tr} F^m\right) \, .
$$

For analytic hyperbolic diffeomorphisms and weights, it has been
known for a long time  that $d_\LL(z)$ is an entire function
when the dynamical foliations are analytic: This is the content of the
fundamental paper of Ruelle \cite{Ru}, who showed that the transfer operators are nuclear on a suitable space of functions  using
Grothendieck's theory  of nuclear operators. More
recently, Rugh \cite{Rugh} and Fried
\cite{Fr} studied $d_\LL(z)$ in the hyperbolic analytic framework, but without any assumption on the
foliations, giving a spectral interpretation of its zeroes (however, not quite as
the eigenvalues of a natural transfer operator $\LL$).
In the case of finite differentiability $r$, the connection
between transfer operators
and dynamical determinants of {\it expanding}
endomorphisms has been well  understood by Ruelle  (see \cite{Ru2}).

In a ground-breaking article \cite{Ki} circulated as a
preprint since 1995, Kitaev considered hyperbolic diffeomorphisms
of finite differentiability $C^r$, and obtained a remarkable formula
$\rho_r:=\inf_{q<0<p, p-q < r-1} \rho^{p,q}(T,g)$ (see Section \ref{2} for a definition 
of $\rho^{p,q}(T,g)$) as a lower bound for the radius of 
a disc in which  $d_\LL(z)$ admits a holomorphic extension. But Kitaev did not construct
a Banach space and his approach does {\it not}
give spectral results. Interpreting the zeroes of $d_\LL(z)$ in the disc furnished by Kitaev
as inverse eigenvalues of a transfer operator remained a challenging problem for over
a decade. 

The main contribution of the present paper is to close this problem
(Theorems~\ref{mainprop} and  ~\ref{main}). Along the way, we give a new
proof of Kitaev's result.
In addition, we give a new variational-like interpretation of Kitaev's formula $\rho^{p,q}(T,g)$
as a kind of generalised topological pressure $Q^{p,q}(T,g)$
(Lemma ~\ref{Kitform}).

\smallskip

{\bf Finding appropriate Banach spaces}

The first reason why this problem remained open for so long is that there 
were until recently no good Banach spaces adapted to the transfer operators of hyperbolic dynamical
systems in finite differentiability:
For a long time, statistical properties of hyperbolic diffeomorphisms were investigated 
using symbolic dynamics via Markov partitions. Since the correspondence is not smoother 
than H\"older, the information thus obtained
on the spectrum of transfer operator was severely limited, and 
this made it difficult to go beyond the results on dynamical zeta functions
by Parry and Pollicott (see \cite{PP}).  (See Section ~2 for a discussion 
about dynamical zeta functions $\zeta_{T,g}(z)$.)
Recently, in a pioneering work\cite{BKL}, Blank, Keller and Liverani 
introduced scales of Banach spaces of distributions
on the manifold adapted to hyperbolic diffeomorphisms and proved that the transfer operators acting on those Banach spaces have a spectral gap. 
However, there were  technical restrictions in the methods in \cite{BKL},
which did not allow them to go beyond Lipschitz smoothness. 
These restrictions were removed by  Gou\"ezel and Liverani \cite{GoL} and 
by the authors \cite{BT} independently, but using different kind of Banach spaces. 
The intuitive idea is the same for both kind of Banach spaces: They consist of distributions
on the manifold, which are as smooth as $C^p$ functions for some $p>0$ in directions close to the
unstable direction, and which are as rough as distributions of order $-q$ for some $q<0$ in directions
close to the stable direction. 
However the real construction of the Banach spaces in \cite{GoL} and \cite{BT} are quite different. We refer to the original papers for details. (The reader-friendly survey \cite{BT2} will be helpful to get ideas in the construction in \cite{BT}.)

In our first main result (Theorem \ref{mainprop}), 
we introduce  yet another scale of Banach spaces, 
$\BB^{p,q}$,
which 
is a kind of hybrid of those in \cite{BT} and \cite{GoL}, and gives a better upper bound 
on the essential spectral radius. This upper bound  coincides with
Kitaev's formula $\rho^{p,q}(T,g)$ (Lemma ~\ref{Kitform}). In view of the results \cite{GL} of Gundlach and Latushkin for expanding maps, we believe that our bound is optimal. 
\smallskip

{\bf Introducing appropriate traces}

The second difficulty to solve this problem in the case of hyperbolic $C^r$ diffeomorphisms is 
to find an appropriate definition for the trace of transfer operators that are not even compact. 
Liverani \cite{Li} found a simple argument to relate 
eigenvalues of $\LL$ with zeroes of the dynamical Fredholm determinant
$d_\LL(z)$, using the Banach spaces in \cite{GoL}. More recently, Liverani and Tsujii \cite{LT} provided an abstract 
argument that is adaptable to both of the Banach spaces in \cite{GoL} and \cite{BT} and that improves the result in \cite{Li}. Still, by technical reasons, the methods 
in \cite{Li} and \cite{LT} give the relation only on a strictly smaller disk (by the factor of one half, at least)  than that given in Kitaev's \cite{Ki} formula. Our second main result (Theorem ~ \ref{main}) improves this point. 

In this paper, we use the structure of our Banach spaces to define the trace. The basic idea in the construction of our Banach spaces is to view functions $u$ on the manifold as superpositions of countably many parts $u_\gamma$, $\gamma\in\Gamma$, each of which is  compactly supported in  
Fourier space (in charts). Accordingly we regard the transfer operator as a countable matrix of operators $\mathcal{L}_{\gamma\gamma'}$. 
Each operator $\mathcal{L}_{\gamma\gamma'}$ turns out to have a smooth kernel. Thus we
may define the trace of $\mathcal{L}_{\gamma\gamma}$ as the integration of its kernel along the diagonal, and then the trace
of the transfer operator $\mathcal{L}$ as the sum of the traces of the $\mathcal{L}_{\gamma\gamma}$.
We found that  hyperbolicity of the diffeomorphism ensures that this trace  
coincides with the expected sum over fixed points. 
Then, using the abstract notion of approximation numbers \cite{Pie}, we estimate the traces thus defined and  get our second main result (Theorem ~ \ref{main}). 
This implement the idea mentioned in the beginning for the case of $C^r$ hyperbolic diffeomorphisms.

\subsection{Main results}
 
In the following, $X$ denotes a $d$-dimensional $C^\infty$ Riemann manifold
and  $T:X \to X$ is a  
diffeomorphism of class $C^r$ for some $r>1$. 
If $r$ is not an integer, 
this means that the derivatives of $T$ of order $[r]$ satisfy a H\"older condition of 
order $r-[r]$. Our standing assumption is that there exists
a hyperbolic basic set $\Lambda\subset X$ for $T$,
that is, a compact $T$-invariant subset that is hyperbolic, isolated and transitive. 
By definition there exist a compact isolating 
neighbourhood $V$ such that $\Lambda=\cap_{m\in \integer}T^{m}(V)$ and 
an invariant decomposition 
$T_{\Lambda}X=E^u\oplus E^s$ of the tangent bundle over $\Lambda$, such that 
$\|DT^m|_{E^s}\|\le C\lambda^{m}$ 
and 
$\|DT^{-m}|_{E^u}\|\le C\lambda^m$, for all $m\ge 0$ and $x\in \Lambda$, with some constants 
$C>0$ and $0<\lambda<1$. 
By transitivity, the dimensions of $E^u(x)$ and $E^s(x)$ are constant, which are denoted 
by $d_u$ and $d_s$ respectively. We will suppose that neither $d_s$ nor $d_u$ is zero.

For $s\ge 0$, let $C^s(V)$ be the set of complex-valued $C^s$ functions on $X$ with support  
contained in the interior of $V$. 
The Ruelle transfer operator with weight $g\in C^{r-1}(V)$ is defined by
\[
\LL=\LL_{T,g}:C^{r-1}(V)\to C^{r-1}(V),\quad \LL \varphi(x)=g(x)\cdot \varphi\circ T(x).
\]
Our first theorem improves the results of \cite{BT} and \cite{GoL} on the
spectrum of $\LL$. 
For a 
$T$-invariant Borel probability measure $\mu$ on $\Lambda$, we write
$h_\mu(T)$ for the metric entropy of $(\mu,T)$, 
and $\chi_\mu(A)\in \real \cup \{-\infty\}$ for the largest Lyapunov
exponent of a linear cocycle $A$ over $T|_\Lambda$,
with $(\log \|A\|)^+\in L^1(d\mu)$.
Let $\MM(\Lambda,T)$ denote the set of  $T$-invariant 
ergodic Borel probability
measures  on  $\Lambda$. Then the theorem is stated as follows.

\begin{theorem}\label{mainprop}
For each real numbers $q<0<p$  so that $p-q<r-1$,
there exists a Banach space $\BB^{p,q}(T,V)$ of distributions on $V$,
containing $C^s(V)$ for any $s>p$, and  contained
in the dual space of $C^{s}(V)$ for any $s>|q|$, with the following
property:

For any $g\in C^{r-1}(V)$, the Ruelle operator $\LL_{T,g}$ 
extends to a bounded operator on $\BB^{p,q}(T,V)$ and the
essential spectral radius of that extension is  not larger than
\begin{align*}
Q^{p,q}(T,g)&=\exp \sup_{\mu \in \MM(\Lambda,T)}
\Bigl \{h_\mu(T) + \chi_\mu\left (\frac{g}{\det (DT|_{E^u})}\right ) \\
\nonumber&\qquad\qquad \qquad\qquad\qquad+
\max\bigl \{p \chi_\mu(DT|_{E^s}), |q| \chi_\mu(DT^{-1}|_{E^u} )\bigr \}
\Bigr \}\, .
\end{align*}
\end{theorem}

Note that, in the setting of $C^r$ expanding endomorphisms, Gundlach and Latushkin \cite[\S 8]{CL},\cite{GL} 
showed that the essential spectral radius of the transfer operator acting on $C^{r-1}(X)$ is 
given exactly by  a variational expression analogous to  $Q^{p,q}(T,g)$.

\begin{remark}\label{lyap}
By upper-semi-continuity of $\mu\mapsto h_\mu$, $\mu\mapsto \chi_\mu(A)$, the supremum 
in the expression for $Q^{p,q}(T,g)$ is a maximum.
Also we have $\chi_\mu(g/\det (DT|_{E^u}))=\int \log |g|\, d\mu
-\int \log |\det( DT|_{E^u})|\, d\mu$. This artificial
expression as a Lyapunov exponent will make sense when we consider Ruelle operators 
on sections of vector bundles in the next section.
\end{remark}

\begin{remark} Note that we have 
\[
 Q^{p,q}(T, g)\le \lambda^{\min\{p,-q\}}\cdot 
Q^{0,0}(T,g)  <Q^{0,0}(T,g)\, .
\] 
We shall see in Remark \ref{re:main} that, if
$g>0$ on $\Lambda$, the spectral radius of
$\LL_{T,g}$ on $\BB^{p,q}(T,V)$  coincides with
$Q^{0,0}(T,g)$.
\end{remark}

\smallskip
To compare the results in this paper with those in Kitaev's  article \cite{Ki}, we 
next give an alternative expression for $Q^{p,q}(T, g)$. 
For $g\in C^{0}(V)$ and $m\ge 0$, we write 
$$ 
g^{(m)}(x)=\prod_{k=0}^{m-1} g(T^k(x)) \, .
$$
We define local hyperbolicity exponents for $x\in \Lambda$ and $m\in \integer_+$
by\footnote{The definition of $\lambda_x(T^m)$ may look a bit strange. We need this formulation for the extension of $E^s(x)$ just below.}
\begin{equation}\label{hypexp2}
\begin{aligned}
\lambda_{x}(T^{m})&=\sup
\left\{ \left.
\frac{\|DT^{m}_x(v)\|}{\|v\|} \;\right |\; \; DT^{m}_x(v)\in E^s(T^m(x))\setminus\{0\}
\right\} \le C\lambda^m, \\
 \nu_x(T^{m})
&=\inf
\left\{\left. 
\frac{\|DT^m_x(v)\|}{\|v\|}\; \right|\;
v\in E^u(x)\setminus\{0\}\right\}\ge C^{-1}\lambda^{-m}.
\end{aligned}
\end{equation}
For real numbers $q$ and $p$, an integer $m\ge 1$ and $x\in \Lambda$, we set
\[
\lambda^{(p,q,m)}(x)=\max\bigl \{ (\lambda_{x} (T^m))^{p}, (\nu_{x} (T^m))^{q}
\bigr \} \, .
\]

We may extend $E^s(x)$ and $E^u(x)$ to continuous bundles on $V$
(which are not  necessarily invariant), so that the inequalities in (\ref{hypexp2}) 
hold for all $x\in \cap_{k=0}^{m-1}T^{-k}(V)$ and for all $m\ge 0$, with some constant $C$. 
Taking such an extension\footnote{The choice of extensions is not essential.}, we extend 
the definition of $\lambda_{x}(T^{m})$, $\nu_x(T^{m})$, and $\lambda^{(p,q,m)}(x)$ to 
$\cap_{k=0}^{m-1}T^{-k}(V)$.
Letting $dx$ denote normalised Lebesgue measure on
$X$, define  for integers $m\ge 1$ and  $p$, $q\in \real$
\begin{equation}\label{2}
\rho^{p,q}(T,g,m)= 
\int_{X}  |g^{(m)}(x)|\cdot  \lambda^{(p,q,m)}(x)\,  dx \, .
\end{equation}
Kitaev\cite{Ki} proved\footnote{Kitaev 
used the notation $\rho^{p,-q}(\LL)$ for our $\rho^{p,q}(T,g)$.} that the limit 
\[
\rho^{p,q}(T, g)=\lim_{m\to \infty}(\rho^{p,q}(T, g, m))^{1/m}
\]
exists for all $q\le 0\le p$ in $\real$ 
 and 
$g\in C^\delta(V)$ with $\delta >0$. In  Section~ \ref{bd}, we show:

\begin{lemma} \label{Kitform}
For any
$g\in C^\delta(V)$ with $\delta>0$, we have
$Q^{p,q}(T,g)= \rho^{p,q}(T, g)$ for all real numbers $q\le 0\le p$.
\end{lemma}

In \cite{BT}  we proved a result  similar to Theorem~\ref{mainprop}, with
$\BB^{p,q}(T,V)$ replaced by other spaces of anisotropic distributions $C^{p,q}_*(T,V)$, respectively
$W_*^{p,q,t}(T,V)$ for $1<t<\infty$,  and 
with the bound $Q^{p,q}(T, g)$ replaced by $R^{p,q,\infty}(T,g)$, respectively
$R^{p,q,t}(T,g)$,  where
\[
R^{p,q,t}(T,g)=\lim_{m\to \infty}\left(
 \sup_{\Lambda} |\det DT^m|^{-1/t} (x)\cdot  |g^{(m)}(x)| \lambda^{(p,q,m)}(x)\right)^{1/m} \, . 
\] 
Note that if $|\det DT|\le 1$ then $\inf_{t\in [1,\infty]} R^{p, q,t}(T,g)= R^{p, q,\infty}(T,g)$. Since 
\begin{align*}
\exp\biggl(\chi_\mu(g) + \max\bigl \{p \chi_\mu(DT|_{E^s}), &|q| \chi_\mu(DT^{-1}|_{E^u} )\bigr \}\biggr)\\
&\qquad  \le 
\lim_{m\to \infty}\left(
 \sup_{\Lambda}  |g^{(m)}(x)| \lambda^{(p,q,m)}(x)\right)^{1/m}\;,
\end{align*}
the variational principle tells that we have $Q^{p,q}(T,g)\le R^{p,q,\infty}(T,g)$ in general and the equality holds only if the supremum in the definition of  
$Q^{p,q}(T,g)$ is attained by the SRB measure for $T$. 
Therefore  Theorem~\ref{mainprop} can be viewed as
an improvement of our previous result \cite{BT}.
In  Appendix~\ref{ineq} we prove that, in general,
\begin{equation}\label{ineqq}
\rho^{p,q}(T,g)\le \inf_{t\in [1,\infty]} R^{p, q,t}(T,g) \,, 
\end{equation}
where the inequality can be strict.

Another improvement on
\cite{BT} is that we now have the same bounds for the 
essential spectral radii of the pull-back operator and the Perron-Frobenius operator, 
which are dual of each other: Take\footnote{We need to  multiply by $h$ 
to localize functions to $V$. If $T$ is Anosov, we may forget about $h$. }
$h\in C^\infty(V)$
so that $h\equiv 1$ on a  neighbourhood of $\Lambda$, and consider the pull-back operator  
$\varphi\mapsto h\cdot \varphi \circ T$ on
$B^{p,q}(T,V)$, and the Perron-Frobenius operator
$\varphi\mapsto (h \cdot \varphi)\circ T^{-1}\cdot |\det (DT^{-1})|$
on $B^{-q,-p}(T^{-1},V)$. Exchanging the
roles of $E^s$ and $E^u$, the bounds 
in Theorem~\ref{mainprop}
for the essential spectral radii of these operators coincide:
\[
Q^{p,q}(T, g)= Q^{-q,-p}(T^{-1}, g\cdot |\det (DT^{-1})|) \, .
\]

\bigskip

We next turn to  dynamical Fredholm determinants. 
The dynamical Fredholm determinant $d_{\LL}(z)$ corresponding to the Ruelle 
transfer operator $\LL=\LL_{T,g}$ is
\begin{equation}\label{eqn:zeta}
d_{\LL}(z)=\exp \left(- \sum_{m=1}^\infty\frac {z^m}{m}
\sum_{T^m(x)=x}
\frac{g^{(m)}(x)} {|\det (\id-DT^m(x))|} \right)\,  .
\end{equation}
The power series in $z$ which is exponentiated
converges only if $|z|$ is sufficiently small. (See Remark \ref{re:main}.) 
Our main result is about the analytic continuation of 
$d_{\LL}(z)$:

\begin{theorem}\label{main}
Let $g\in C^{r-1}(V)$.

\noindent {\rm (1)} The function $d_{\LL}(z)$ extends  holomorphically
to the disc of radius $(Q_{r-1}(T,g))^{-1}$ with
\[
Q_{r-1}(T,g)=\inf_{q<0<p, \, p-q<r-1} Q^{p,q}(T,g)\, .
\]
{\rm (2)} For any real numbers $q<0<p$  so that $p-q<r-1$,
and each $z$ with $|z|< (Q^{p,q}(T, g))^{-1}$, 
we have $d_\LL(z)=0$  
if and only if $1/z$ is an  eigenvalue of $\LL$ on $\BB^{p,q}(T,V)$, 
and the order
of the zero coincides with the algebraic multiplicity of the eigenvalue.
\end{theorem}

\begin{remark}\label{re:main}
The  sum over $m$ in the right hand side of (\ref{eqn:zeta}) converges when 
\[
|z|< 
\exp
\biggl (
-P_{top}\bigl (T|_\Lambda,\log (|g|/|\det( DT|_{E^u})| )
\bigr )
 \biggr )
=(Q^{0,0}(T,g))^{-1}\, ,
\]
so that $d_{\LL}(z)$ is a nowhere vanishing holomorphic function in this disc. 
To see this, note that there is $C\ge 1 $ so that for all $m$ and all
$x\in \Lambda$ with $T^m(x)=x$
$$ 
C^{-1}\le \frac{|\det (\id-DT^m(x))|}{ |\det( DT^m|_{E^u})(x)|} \le C\, ,
$$
then use the Cauchy  criterion for the convergence
of a power series and the expression of topological pressure as an asymptotic
weighted sum over periodic orbits (see, e.g., \cite[Prop. 5.1]{PP}).
If $g> 0$ on $\Lambda$, then it follows from Pringsheim's theorem
on power series with positive coefficients \cite[\S17]{La}
that $d_{\LL}(z)$ has a
zero at $(Q^{0,0}(T,g))^{-1}$. 
\end{remark}

\medskip

This paper is organized as follows. 
In Section~\ref{add}, we discuss about transfer operators acting
on sections of vector bundles, with applications to dynamical zeta
functions.
In Section~\ref{bd}, we present a
key alternative expression for the bound 
$Q^{p,q}(T,g)$ (useful also to prove
both main theorems), and we
prove Lemma~\ref{Kitform}.

In Section~\ref{sec:spaces}, we consider the transfer operator $L$ on $\real^d$ for  a $C^r$ diffeomorphism $\TT$ and a $C^{r-1}$ weight $G$. 
We first introduce the Banach space $\BB^{\Theta,p,q}(K)$ of anisotropic distributions on a compact subset $K\subset \real^d$, slightly modifying the definitions in \cite{BT}: the $L^\infty$ norm in the definition
of anisotropic spaces in \cite{BT} is replaced by a 
mixed norm, which involves both the supremum norm and
the $L^1$-norm along manifolds close to unstable manifolds.
To study the action of the transfer operator $L$ on this Banach space, we work with an auxiliary operator $M$, which is an infinite matrix of 
operators describing transitions induced by
$L$ between frequency bands in Fourier
space. We observe that the operator $M$ is naturally decomposed as $M_b + M_c$
 with $M_b$ having small spectral radius and $M_c$ a compact operator.
In Lemma~\ref{lm:Mb}, we give a simple estimate on the operator norm of $M_b$. In Subsection~\ref{appr}, we study the approximation numbers of $M_c$  and show, in particular, that $M_c$ is compact. 
The use of approximation numbers to study dynamical transfer operators seems to be new.

In Section~\ref{sec:L}, we introduce
the anisotropic Banach spaces $\BB^{p,q}(T,V)$, and prove Theorem \ref{mainprop}. 
Taking a system of local charts on $V$ adapted to hyperbolic structure of $T$, we consider the system $\KK$ of transfer operators that $\LL$ induces on the local charts. Then we associate an auxiliary operator\footnote{It is possible to work directly with $\LL$, decomposing it
into a compact term $\LL_c$ and a bounded term  $\LL_b$, on $\BB^{p,q}(T,V)$, 
in the spirit of \cite{BT2}. Then the flat trace of $(\LL^m)_b$ is not zero, but
it decays exponentially,  arbitrarily fast \cite{Bbook}.}$\MM$ to $\KK$, in the same manner as we associate $M$ to $L$ in Section \ref{sec:spaces}. 
The spectral data of $\KK$ and $\MM$ turn out to be (almost) identical with that of $\LL$. 
We will decompose $\MM^m$ for $m\ge 1$ as $(\MM^m)_b+(\MM^m)_c$,
where $(\MM^m)_c$ is compact and $(\MM^m)_b$ has norm
smaller than $C(Q^{p,q}(T,g)+\epsilon)^m$, proving Theorem \ref{mainprop}.  

In Section \ref{sec:trace}, we introduce a formal trace $\tr^\flat(\PP)$,  called the flat trace,
and a formal determinant 
$\det^\flat(\id-z\PP)=\exp-\sum_{m\ge 1}\frac{z^m}{m}\tr^\flat(\PP^m)$.
The flat trace is
a key tool inspired from \cite{BB,BR}.
(The terminology was borrowed from 
 Atiyah--Bott \cite{AB2}, but we do not
relate our flat trace to theirs.)
Our flat trace coincides, on the one hand, with the usual trace for finite rank operators and, on the other hand, with the dynamical trace for each $\MM^m$:
$$
{\tr}^\flat(\MM^m)
=\sum_{T^mx=x}\frac{g^{(m)}(x)}
{|\det (\id-DT^m(x))|}\, ,\quad
\mbox{so } \quad d_\LL(z)={\det}^\flat(\id-z\MM)\,.
$$ 
Also, the flat trace $\tr ^\flat((\MM^m)_b)$ vanishes for all large enough $m$.

In Section \ref{S5}, we give the proof of Theorem \ref{main}.
The basic idea of the proof is then to exploit the formal determinant identity\footnote{The operator $\DD(z)=z\MM_c(\id-z\MM_b)^{-1}$ can be viewed
as a kneading operator, \cite{BR}, \cite{Bai}.}
\begin{equation}\label{knead}
{\det}^\flat(\id-z\MM)={\det}^\flat(\id-z\MM_c(\id-z\MM_b)^{-1})\cdot {\det}^\flat(\id-z\MM_b) \, .
\end{equation}
If $r>d+1+p-q$ each operator $(\MM^m)_c$ 
turns out to be  an operator with summable approximation numbers, and 
our proof in this case is fairly simple, although we cannot apply
(\ref{knead}) directly, since we only know that $\tr ^\flat((\MM^m)_b)=0$
and that the spectral radius of $(\MM^m)_b$
is smaller than $(Q^{p,q}(T,g)+2\epsilon)^m$ {\em for large $m$.} 
If $r\le d+1+p-q$, we need more estimates 
since only some iterate 
of $(\MM^m)_c$ has summable approximation
numbers. Still the proof is straightforward. 
\footnote{See \cite{Bbook} for a ``regularised determinant" alternative
to the argument in Section 7.}

In Appendix \ref{apd:eigenvalue}, we discuss about eigenvalues and eigenvectors of the transfer operator $\LL$ on  different Banach spaces.


\section{Operators on vector bundles and dynamical zeta functions}
\label{add}

We may generalize the statements and proofs of the main results 
to similar operators acting on spaces of sections of vector bundles. 
Since  Ruelle zeta function is given as a product of the dynamical 
Fredholm determinants of such operators \cite{Fr0, Ru}, we can derive statements for Ruelle zeta functins from our main theorems.  
See also \cite{PP} for a presentation of classical results about
dynamical zeta functions.

For $r>1$, $T$, and $V$ as in Section ~1,
let $\pi_{B}:B\to V$ be a finite dimensional complex vector bundle, and let 
$\TTT:B\to B$ be a $C^{r-1}$ vector bundle endomorphism such that 
$\pi_{B}\circ \TTT=T^{-1}\circ \pi_B$. Denote the natural action of $\TTT$ on continuous 
sections 
of $B$ by $\LL=\LL_{\TTT}$, that is, 
$\LL u(x)=\TTT(u(T(x)))$.  
Then we can define $Q^{p,q}(T,\TTT)$ in parallel with the definition of $Q^{p,q}(T,g)$ in 
Section~ 1, replacing  $\chi_\mu(g/\det(DT|E^u))$ by $\chi_\mu(\TTT/\det(DT|E^u))$. 
Putting, for $m\ge 1$,
\[
|\TTT^{(m)}|(x)=\|\TTT^m_x:B_x\to B_{T^{-m}(x)}\| \, ,
\]
we can define $\rho^{p,q}(T,\TTT,m)$  by using the same formal expression as for
$\rho^{p,q}(T,g,m)$. 

The next statement is just a formal extension of Theorem \ref{mainprop} and
Lemma~\ref{Kitform}:

\begin{theorem}\label{mainpropV}
Let $q<0<p$ be so that $p-q<r-1$.
There exists a Banach space $\BB^{p,q}(T,B)$ of distributional sections of $B$,
containing $C^s$ sections for any $s>p$, so that 
the operator $\LL_\TTT$ 
extends to a bounded operator on $\BB^{p,q}(T,B)$, and its
essential spectral radius on this space is 
not larger than $Q^{p,q}(T, \TTT)=\rho^{p,q}(T,\TTT)$.
\end{theorem}
Note that if $B$ is the $k$-th exterior power of the
cotangent bundle of $X$ then  $\BB^{p,q}(T,B)$ is a space of
 currents on $X$.
 
The dynamical Fredholm determinant of $\LL=\LL_{\TTT}$ as above is defined by
\[
d_{\LL}(z)=\exp - \sum_{m=1}^\infty\frac {z^m}{m}
\sum_{T^m(x)=x}
\frac{\tr \,\TTT_x^m} {|\det (\id-DT^m(x))|} \, .
\]
A formal extension of Theorem \ref{main} gives:
\begin{theorem}\label{mainV}
For any $q<0<p$  so that $p-q<r-1$, 
the function $d_{\LL}(z)$ extends  holomorphically
to the disc of radius $(Q^{p,q}(T,\TTT))^{-1}$, and its zeroes in 
this disc are
exactly  the inverses of the eigenvalues of $\LL_{\TTT}$ on $\BB^{p,q}(T,B)$, 
the order
of the zero coinciding with the multiplicity of the eigenvalue.
\end{theorem}

Let $\pi_L:L\to \Lambda$ be the orientation line bundle for the  bundle 
$\pi_{E^u}:E^u\to \Lambda$, that is, the fiber of $L$ at $x\in B$ is 
isomorphic to the real line whose unit vectors corresponding to an orientation on $E^u(x)$. 
By shrinking the isolating neighbourhood $V$, we may extend it to a continuous line bundle $\pi_L:L\to V$. 
Let $g\in C^{r-1}(V)$. 
For $k=0, 1,\cdots, d$, let $\pi:B_k=(\wedge^k T^*X)\otimes L\to V$ and let $\TTT_{k}:B_k \to B_k$ be the vector bundle endomorphism defined by $
\TTT_k(w)=(g\circ \pi)\cdot T^*(w)$. Let $\LL_k$ be the natural action of $\TTT_k$ on the sections of $B_k$. 
Then the Ruelle zeta function 
\begin{equation}
\zeta_{T,g}(z)=\exp\left(  \sum_{m=1}^\infty\frac {z^m}{m}
\sum_{T^m(x)=x}g^{(m)}(x) \right)\, .
\end{equation}
can be written as 
\[
\zeta_{T,g}(z)=\prod_{k=0}^{d} d_{\LL_k}(z)^{(-1)^{k+\dim E^u+1}} \, .
\]
Thus we obtain as a corollary of Theorem~\ref{mainV}:

\begin{corollary} The Ruelle zeta function $\zeta_{T,g}(z)$ extends as a meromorphic function 
to the disk of radius
\[
\min_{0\le k\le d}\;\sup\;\left\{
Q^{p,q}(T,\TTT_k)^{-1} \;|\; q<0<p,\; p-q<r-1\right\}\, .
\]
\end{corollary}



\section{Alternative expressions for the bound $Q^{p,q}(T,g)$}
\label{bd}

In this section, we introduce two more expressions, $Q^{p,q}_*(T,g)$ and $\rho^{p,q}_*(T,g)$, in addition to $Q^{p,q}(T,g)$ and $\rho^{p,q}(T,g)$, inspired by \cite{Ki}. 
And we show that these four expressions are all equivalent, proving Lemma~\ref{Kitform} especially. Along the way, we express $\log Q^{p,q}(T,g)$ 
as  a double limit of  topological pressures 
(Lemma~\ref{QQ_0}). Note that the expression  $Q^{p,q}_{*}(T,g)$ will play a central role in the proofs of Theorems~ \ref{mainprop}
and ~\ref{main}
in the following sections.

In this section, $r>1$, $T$ and $\Lambda\subset V$
are as in Section~1, but  we only assume 
$g\in C^\delta( V)$ for some $\delta>0$ (and sometimes
only that $g\in C^0(X)$), even 
if $r$ is large. 
\begin{remark}
Unlike the standard argument\cite{W} on topological pressure, we consider the case where the function $g$ may vanish at some points on $\Lambda$. If we assumed that $g$ vanishes nowhere on $\Lambda$, the argument in this section should be simpler and partly follow form the standard argument. 
\end{remark}

\subsection{The expression $Q^{p,q}_*(T,g)$ and topological pressure}
\label{3.2}
Recall that, in Section~ \ref{sec:intro}, we extended the decomposition 
$T_xX=E^s(x)\oplus E^u(x)$ on $\Lambda$ to $V$ and defined $\lambda_x(T^m)$ and $\nu_x(T^m)$ 
for $x\in \cap_{k=0}^{m} T^{-k}(V)$. 
Using this extension, we also define 
\begin{equation}\label{defdet}
|\det (DT^m|_{E^u})|(x)\, \quad \mbox{ for } x\in \cap_{k=0}^{m} T^{-k}(V)\, ,
\end{equation} 
as 
the  expansion factor of the linear mapping $DT^m:E^u(x)\to DT^m(E^u(x))$, with respect to 
the volume induced by the Riemannian metric on each $d_u$-dimensional linear subspace. 
Note that, for each $g\in C^0(V)$, the sequences of functions
$g^{(m)}$ and $|\det (DT^m|_{E^u})|$ are multiplicative, while 
$\lambda^{(p,q,m)}$ is submultiplicative in $m$ for all
real numbers $q\le 0\le p$. 
In particular,
$|g^{(m)}| \cdot \lambda^{(p,q,m)} \cdot |\det (DT^m|_{E^u})|^{-1}$ 
is submultiplicative
in $m$ for such $p$ and $q$. 


We say that $\WW$ is a cover of $V$ if it is a finite
cover $\WW=\{W_i\}_{i \in \II}$ of $V$
by open subsets of $X$ and if, in addition, 
the union $\cup_{i\in \II}W_i$ is contained in a compact isolating neighbourhood $V'$ of $\Lambda$.
For such a cover $\WW$ and integers $n< m$, put
\[
\WW^m_n=\{\cap_{k=n}^{m-1} T^{-k}(W_{i_k})\mid (i_k)_{k=n}^{m-1} \in \II^{m-n}  \} \, ,
\]
and set $\WW^m=\WW^m_0$ for $m\ge 1$. 
Then $\WW^m$ is a cover of $V^m:=\cap_{k=0}^{m-1}T^{-k}(V)$.
We say that a cover $\WW$ of $V$
is generating  if  the diameter of  $\WW^m_{-m}$ 
tends to zero as $m\to \infty$. (Generating covers exist
because $\cap_{k=-m}^m T^{-k} V$ is contained in a small
neighbourhood of $\Lambda$ for large $m$.)
For real numbers
$p$ and $q$, an integer $m\ge 1$, a generating cover $\WW$
of $V$, and $g \in C^0(X)$, we define
\begin{align}\label{defQ_0}
Q^{p,q}_{*}(T,g,\WW,m)&=\min_{\WW'}
\left(\sum_{U\in \WW'}
\sup_{U} \frac{|g^{(m)}| \lambda^{(p,q,m)}}
{|\det (DT^m|_{E^u})|} \right)
\end{align}
where the minimum $\min_{\WW'}$
 is taken over subcovers $\WW'\subset \WW^m$ of $V^m$. 
By sub-multiplicativity with respect to $m$,  the following limits exist
if $q\le 0\le p$:
\begin{align*}
Q^{p,q}_{*}(T,g,\WW)&=\lim_{m\to \infty} 
\left( Q_{*}^{p,q}(T,g,\WW,m)\right)^{1/m}\, .
\end{align*}

The following lemma may not be new. But, since we did not find it in the literature, we provide 
a proof.

\begin{lemma}\label{note}
For any generating cover $\WW$ of $V$
and 
$g\in C^0(X)$ with $\inf_X |g|>0$, we have
$
\log Q^{0,0}_{*}(T,g,\WW)
=P_{top}\bigl (T|_\Lambda, \log(|g|/|\det (DT|_{E^u})|) \bigr )
$.
\end{lemma}
\begin{proof}
It is enough to show 
\begin{equation}\label{desired}
\log Q^{0,0}_*(T,g,\WW)\le P_{top}(T|_{\Lambda},\log (|g|/|\det(DT|_{E^u})|))\, ,
\end{equation}
since the inequality in the opposite direction is clear. 
Let $\WW=\{W_{i}\}_{i\in\II}$. Take another cover  
$\UU=\{U_{i}\}_{i \in \II}$ of $V$,
so that $U_{i}\Subset W_{i}$ for $i\in \II$. Consider small $\epsilon>0$ so that, 
for each $i\in \II$, the $\epsilon$-neighbourhood of $U_{i}$ is contained in 
$W_i$.

Let $W_{\vec \imath}:=\bigcap_{k=0}^{m-1}T^{-k}(W_{i_k})$ and 
$U_{\vec \imath}:=\bigcap_{k=0}^{m-1}T^{-k}(U_{i_k})$ 
for $\vec\imath=(i_k)_{k=0}^{m-1}\in \II^m$. 
For each $m\ge 1$, let $Q_\Lambda(T,g,\UU,m)$ be the minimum of 
\[
\sum_{\vec \imath\in \II'} 
\sup_{U_{\vec \imath}\cap \Lambda} \frac{|g^{(m)}|}
{|\det(DT^m|_{E^u})|} 
\]
over subsets $\II' \subset \II^m$ such that $
\{U_{\vec \imath}\cap \Lambda 
\mid {\vec \imath} \in \II' \}$ is a cover of $\Lambda$. 
Let $\II'=\II'(m)$ be a subset of $\II^m$ that attains this minimum. 

Since $V$ is 
an isolating neighbourhood for the hyperbolic basic set $\Lambda$, we can take large $N$ so that, if $T^{k}(x)\in V$ for 
$0\le k\le n+2N$, there exists a point $y\in \Lambda$ such that $d(T^{N+k}(x), T^k(y))<\epsilon$ for all $0\le k< n$. This implies that
\[
\{W_{\vec \imath}\mid  \vec \imath=(i_k)_{k=0}^{m+2N-1}\in \II^{m+2N}\mbox{ and }  
(i_{k+N})_{k=0}^{m-1}\in \II'(m)\}\subset \WW^{m+2N}
\]
is a cover of $V^{m+2N}$. 
Therefore we have, for all $m\ge 0$,  
\[
\min_{\WW'\subset \WW^{m+2N}}
\left(\sum_{U\in \WW'}
\inf_{U} \frac{|g^{(m+2N)}|}
{|\det (DT^{m+2N}|_{E^u})|} \right)\le C\cdot Q_\Lambda(T,g,\UU,m)\, ,
\]
where the minimum 
 is taken over subcovers $\WW'\subset \WW^{m+2N}$ of $V^{m+2N}$, and hence
\[
\varlimsup_{m\to \infty}\frac{1}{m}\log \min_{\WW'\subset \WW^{m}}
\left(\sum_{U\in \WW'}
\inf_{U} \frac{|g^{(m)}|}
{|\det (DT^{m}|_{E^u})|} \right)\le
P_{top}\left(T|_{\Lambda},\log 
\frac{|g|}{|\det(DT|_{E^u})|}\right).
 \]
Since $g$ is continuous and positive and since $\WW$ is a generating cover, the left hand side coincides with $\log Q^{0,0}_*(T,g,\WW)$. \end{proof}

We next express $\log Q_{*}^{p,q}(T,g,\WW)$ as a limit
of topological pressures
under the  condition $\inf_X |g|>0$:
\begin{lemma}\label{lm:tp}  
If
$\WW$ is a  generating cover of $V$ and if $q\le 0\le p$, then
for  each $g\in C^0(X)$ such that $\inf_X |g|>0$, we have
\begin{equation}\label{top1}
 \log {Q}^{p,q}_{*}(T,g,\WW)
=\lim_{m\to \infty} \frac{1}{m}
P_{top}\biggl (
T^m|_\Lambda, \log 
\frac{|g^{(m)}|\cdot \lambda^{(p,q,m)}}{ |\det (DT^m|_{E^u})|}
\biggr ) \, .
\end{equation}
\end{lemma}

\begin{proof}
The topological pressures in the claim are well-defined because for each $m$
the function $\log h_m$, with
\[
h_{m}:=|g^{(m)}|\cdot  \lambda^{(p,q,m)} \cdot |\det (DT^m|_{E^u})|^{-1}\, ,
\]
is continuous on $\Lambda$.  
The limit in (\ref{top1}) exists by sub-multiplicativity of $m\mapsto h_m$.

For each $\epsilon>0$, there exists $m\ge 1$ so that 
\[
(Q^{p,q}_{*}(T,g,\WW)+\epsilon)^m \ge  
Q^{p,q}_{*}(T,g,\WW,m)= Q^{0,0}_*(T^m,|g^{(m)}| \lambda^{(p,q,m)},\WW^m,1) \, .
\]
By Lemma~\ref{note}, the right-hand side is not smaller than 
$\exp(P_{top}(T^m|_\Lambda, \log h_m))$. Hence 
\[
{Q}^{p,q}_{*}(T,g,\WW)\ge 
\lim_{m\to \infty} \exp((1/m)P_{top}(T^m|_\Lambda, \log h_m))-\epsilon \, .
\]
Since $\epsilon>0$ is arbitrary, this give the inequality in one direction. 

We next show the inequality in the opposite direction. 
By sub-multiplicativity and Lemma~\ref{note}, we have, for any integer $m>0$, that 
\begin{align*}
\log {Q}^{p,q}_{*}(T,g,\WW)&= \lim_{k\to \infty}
\frac{1}{mk}\log Q^{0,0}_*(T^{mk},|g^{(mk)}| \lambda^{(p,q,mk)},\WW^{mk},1)\\
&\le \lim_{k\to \infty}
\frac{1}{mk}\log Q^{0,0}_*(T^{m},|g^{(m)}| \lambda^{(p,q,m)},\WW^{m},k)\\
&=\frac{1}{m}\log Q^{0,0}_*(T^{m}, |g^{(m)}| \lambda^{(p,q,m)},\WW^{m}) =\frac{1}{m}P_{top}(T^{m}|_{\Lambda},\log h_{m})\, .
\end{align*}
This gives the inequality in the opposite direction. 
\end{proof}

To get rid of the assumption $\inf_X |g|>0$, we shall use the following:
\begin{lemma}\label{lm:liminf}
Let $\WW$ be a generating cover of $V$,  let $g\in C^0(X)$, and let $q\le 0\le p$. 
If $g_n$ is a  sequence of 
functions in  $C^0(X)$ so that 
$\inf_X g_n>0$ with  $g_{n}\ge g_{n+1}\ge |g|$ for all
$n$, and $\lim_{n\to\infty}\|g_n- |g|\|_{L^\infty(V)}=0$,
 then 
\[
\lim_{n\to \infty}Q^{p,q}_{*}(T,g_n,\WW)=Q^{p,q}_{*}(T,g,\WW)\, .
\]
\end{lemma}

By Lemmas~ \ref{lm:tp} and~ \ref{lm:liminf}, the exponent
$Q^{p,q}_{*}(T,g,\WW)$ for any $g\in C^0(X)$ does not depend on the generating cover $\WW$. 
So it will be denoted by $Q^{p,q}_*(T,g)$.

\begin{proof} We have only to show 
$\lim_{n\to \infty}Q^{p,q}_*(T,g_n,\WW)\le Q^{p,q}_*(T,g,\WW)$. 
For any $\epsilon>0$, we take  large $m$ such that 
$
Q^{p,q}_*(T,g,\WW,m)\le (Q^{p,q}_*(T,g,\WW)+\epsilon)^m$.
Then take $n_0$ such that $Q^{p,q}_*(T,g_n,\WW,m)\le (Q^{p,q}_*(T,g,\WW)+2\epsilon)^m$ 
for $n\ge n_0$. 
By  sub-multiplicativity, we get
$
Q^{p,q}_*(T,g_n,\WW)\le Q^{p,q}_*(T,g,\WW)+2\epsilon 
$
for  $n\ge n_0$.  
\end{proof}

\subsection{A variational principle}
\label{3.2'}

Lemmas~ \ref{lm:tp} and \ref{lm:liminf} allow us to prove:

\begin{lemma}\label{QQ_0}
$Q^{p,q}(T,g)=Q^{p,q}_*(T,g)$ for $q\le 0\le p$ and $g\in C^0(X)$.
In particular, for every sequence $g_n$ as in Lemma~\ref{lm:liminf}, we have
\begin{equation*}
\log Q^{p,q}(T,g)=\lim_{n\to \infty}\lim_{m \to \infty}
\frac{1}{m}
P_{top}\bigl (
T^m|_\Lambda, \log (g_n^{(m)} \lambda^{(p,q,m)} |\det (DT^m|_{E^u})|^{-1})
\bigr )\, .
\end{equation*}
\end{lemma}
\begin{proof}
We first show the claim when $\inf_X |g|>0$.
For simplicity, we put 
\[
\chi_p(\mu)=p \cdot \chi_\mu(DT|_{E^s}),\qquad 
\chi_q(\mu)=|q|\cdot \chi_\mu(DT^{-1}|_{E^u})
\]
and 
\[
P(\mu)=h_\mu(T) + \int \frac{|g|}{\det (DT|_{E^u})} \, d\mu +
\max \{\chi_p(\mu),  \chi_q (\mu)\}\, , 
\]
so that
$\log Q^{p,q}(T,g)= \sup_{\mu\in \MM(\Lambda,T)} P(\mu)$. 
Next we put
\[
P_m(\mu)=m h_\mu(T) +\int 
\log \left(|g^{(m)}|\cdot \lambda^{(p,q,m)} \cdot |\det(DT^m|_{E^u})|^{-1}\right) \, d\mu \, .
\]
By the variational principle \cite{W}, Lemma \ref{lm:tp}
implies 
\[
\log Q^{p,q}_*(T,g)=\lim_{m\to \infty} 
\sup_{\mu\in \mathcal{M}(\Lambda,T)} \frac{1}{m}P_m(\mu) \, .
\]
Note that, for any invariant probability measure $\mu$, Oseledec's theorem\cite{W} gives
\begin{equation}\label{conv}
\lim_{m\to \infty} 
\frac{1}{m}\log \lambda^{(p,q,m)}(x) =\max \{\chi_p(\mu), \chi_q(\mu) \}
\qquad\mbox{for $\mu$-a.e. $x$.}
\end{equation}

We first show $Q^{p,q}(T,g)\le Q^{p,q}_*(T,g)$. 
There exists a measure $\mu_0\in \mathcal{M}(\Lambda,T)$ such that 
$P(\mu_0)=\log (Q^{p,q}(T,g))$. (See Remark \ref{lyap}.) By (\ref{conv}), we obtain
\[
\log Q^{p,q}_*(T,g)\ge  \lim_{m\to \infty}\frac{1}{m}P_m(\mu_0)=P(\mu_0)=\log Q^{p,q}(T,g) \, .
\]
We next show $Q^{p,q}(T,g)\ge Q^{p,q}_*(T,g)$. For each $m$, we take  
$\mu_m\in \mathcal{M}(\Lambda,T)$ such that 
$P_m(\mu_m)=\sup_{\mu\in \mathcal{M}(\Lambda,T)} P_m(\mu)$.
Then we take  a subsequence $m(i)\to \infty$ such that $\mu_{m(i)}$ converges weakly 
to an invariant probability  measure $\mu_\infty$ on $\Lambda$. 
By decomposing $\mu_{\infty}$ into ergodic components 
we see that  $\log Q^{p,q}(T,g)\ge  P(\mu_\infty)$. Thus, 
the claim $Q^{p,q}(T,g)\ge Q^{p,q}_*(T,g)$ follows if we  show
\begin{equation}\label{eqn:Pm}
P(\mu_{\infty})\ge \lim_{i\to \infty}\frac{1}{m(i)}P_{m(i)}(\mu_{m(i)})\, .
\end{equation}
By the upper 
semi-continuity of entropy, we have 
$h_{\mu_\infty}(T)\ge\lim_{i\to \infty}h_{\mu_{m(i)}}(T)$. 
By sub-multiplicativity of $\lambda^{(p,q,m)}$ and (\ref{conv}), we have
\begin{align*}
\limsup_{i\to \infty}
\int \frac{\log \lambda^{(p,q,m(i))}}{m(i)} d\mu_{m(i)}
&\le 
\inf_{m\ge 1}\int \frac{\log \lambda^{(p,q,m)}}{m} d\mu_{\infty} \le \max\{\chi_p(\mu_\infty), \chi_q(\mu_\infty)\}.
\end{align*}
Therefore we get the inequality (\ref{eqn:Pm}). 

Finally we consider the case $\inf_X |g|=0$. Take
a sequence  $g_n$ as in Lemma~\ref{lm:liminf}.
In view of Lemma~\ref{lm:liminf} and the argument above, it remains to show 
$Q^{p,q}(T,g)=\lim_{n\to \infty}Q^{p,q}(T,g_n)$   
for 
$q\le 0\le p$. Note that the sequence $Q^{p,q}(T,g_n)$ is decreasing and we have 
$Q^{p,q}(T,g)\le \lim_{n\to \infty}Q^{p,q}(T,g_n)$ obviously. 
We show the inequality in 
the opposite direction.
We write $P(g,\mu)$ for $P(\mu)$.
For each $n$, take $\mu_n\in \mathcal{M}(\Lambda,T)$ such that 
$P(g_n,\mu_n)=Q^{p,q}(T,g_n)$ and then take a subsequence $n(i)\to \infty$ so that 
$\mu_{n(i)}$ converges weakly to some invariant probability measure $\mu_\infty$ on $\Lambda$. Then, by 
 upper-semi-continuity
of the entropy and of the largest Lyapunov exponent as a function of
$\mu$, we obtain 
\[
\lim_{n\to \infty} Q^{p,q}(T,g_n)\le 
\liminf_{n\to\infty} P(g_n,\mu_n)\le P(g,\mu_{\infty})\le Q^{p,q}(T,g) \, .
\]
\end{proof}

We may now complete the first step towards the proof of 
Lemma~\ref{Kitform}:

\begin{lemma}\label{step1}
$\rho^{p,q}(T,g)
\le Q_*^{p,q}(T,g)$ for any  $q\le 0\le p$ and $g \in C^0(V)$.
\end{lemma}

\begin{proof}
Take a generating
cover $\WW=\{W_i\}$ of $V$. 
Then, by a standard argument on hyperbolicity, we can show\footnote{To see this, we can use the 
``pinning coordinates" in \cite[\S3.3 and p. 163]{Ki}.} that 
the Riemann volume of $U\in \WW^m$ is bounded by 
$C/|\det(DT^m|_{E^u})(x)|$ for any $x\in U$, where $C$ is a constant 
that does not depend on $U$, $x$ or $m$. 
Then we have, for any subcover $\WW'\subset \WW^m$ of $V^m$,
\[
\rho^{p,q}(T,g,m)\le 
\sum_{U \in \WW'}
\int_{U} 
|g^{(m)}(x)| \lambda^{(p,q,m)}(x)\, dx\le C\cdot \sum_{U \in \WW'}
\sup_{y\in U} 
\frac{|g^{(m)}(y)| \lambda^{(p,q,m)}(y)}{|\det(DT^m|_{E^u})(y)|}.
\]
This implies $\rho^{p,q}(T,g,m)
\le C Q^{p,q}_*(T,g,\WW ,m)$ and hence the lemma.
\end{proof}


\subsection{The expression $\rho^{p,q}_*(T,g)$} 
\label{defpart}
We next introduce an exponent $\rho^{p,q}_*(T,g)$  due
to Kitaev \cite{Ki},
using  partitions of unity. 
A finite family
$\Phi=\{\phi_\omega\}_{\omega\in \Omega}$ 
of $C^\infty$ functions on $X$ is called a partition of unity  for $V$ if 
$0\le \phi_\omega(x)\le 1$ on $X$,  and 
$\sum_{\omega\in \Omega}\phi_\omega(x)\equiv 1$ on $V$. 
The diameter of $\Phi$ is 
$\max_{\omega\in \Omega} \{\mbox{diam}\, \bigl (\supp(\phi_\omega) \cap V \bigr )\}$.
For a partition of unity $\Phi$  and an integer $m\ge 1$, set
\[
\Phi^m=\biggl\{\,\biggl.  \prod_{k=0}^{m-1}\phi_{\omega_k}(T^k(x))\;\biggr| \; 
(\omega_k)_{k=0}^{m-1}\in \Omega^m\biggr\} \, ,
\]
which is a partition of unity for $\cap_{k=0}^{m-1}T^{-k}(V)$. 
For $g\in C^0(V)$, the sequence
\[
\rho^{p,q}_*(T,g,\Phi,m)=\sum_{\phi\in \Phi^m} \|\phi\cdot g^{(m)}\cdot
\lambda^{(p,q,m)}\cdot \det (DT^m| _{E^u})^{-1}\|_{L^\infty}\, ,
\]
is then submultiplicative with respect to $m$ if
$q\le 0\le p$, so that we may put
\[
\rho^{p,q}_*(T,g,\Phi)=\lim_{m\to \infty} \left(\rho^{p,q}_*(T,g,\Phi,m)\right)^{1/m} \, .
\]
An important estimate due to  Kitaev is:

\begin{lemma}[Kitaev \cite{Ki}]\label{lm:Ki12}
Let $q\le 0\le p$.
For  every   partition of unity $\Phi$ for $V$ of sufficiently small diameter
and each $g\in C^\delta(V)$ with $\delta>0$, we have
\begin{equation}\label{rhoKi}
\rho^{p,q}_*(T,g,\Phi)=\rho^{p,q}(T,g)\, .
\end{equation}
\end{lemma}
This lemma implies that $\rho^{p,q}_*(T,g,\Phi)$ takes a constant value for any sufficiently fine  partition of unity $\Phi$. This value is denoted by $\rho^{p,q}_*(T,g)$.
\begin{remark}
In \cite[Lemma 2]{Ki}, the corresponding claim is actually stated for "regular mixed transfer operator (MTO)". 
To get Lemma \ref{lm:Ki12}, we apply that 
claim to the regular MTO induced by $T$ and $g$, using local charts and partitions of unity. See \cite{Ki} and Remark \ref{MTO}. 
\end{remark}

\begin{remark}\label{re:Ki12}
In Lemma \ref{lm:Ki12}, we can prove $\rho^{p,q}_*(T,g,\Phi)\ge \rho^{p,q}(T,g)$  without much difficulty, using the argument as in the proof of Lemma \ref{step1}. But  the inequality in the opposite direction and exactness of the limit in the definition of $\rho^{p,q}(T,g)$ are not easy to prove. 
In general, the functions $\lambda^{(p,q,m)}(x)$ for large $m$ depend on $x$ irregularly,
so that we may not use a simple argument. 

\end{remark}

We finally prove Lemma~\ref{Kitform}:

\begin{proof}[Proof of Lemma~\ref{Kitform}]
By Lemma~\ref{QQ_0} and~\ref{step1},
we have only to show that $\rho^{p,q}(T,g) \ge Q^{p,q}_*(T,g)$. 
We start by a preliminary observation: For any integer $k\ge 1$, we have
\begin{equation}\label{mul}
Q^{p,q}_*(T^k,g^{(k)})=(Q^{p,q}_*(T,g))^k
\quad \mbox{and}\quad
\rho^{p,q}(T^k, g^{(k)})= (\rho^{p,q}(T,g))^k\, .
\end{equation}
The former follows from Lemma~\ref{QQ_0}. 
The latter is a consequence of the definition.

We take a partition of unity 
$\Phi=\{\phi_\omega\}_{\omega\in \Omega}$ of small diameter
so that the intersection multiplicity of the supports of $\phi_\omega$ is less 
than some constant $N_d$ that depends only on the dimension $d$ of $X$.
Then  $\WW=\{\phi_{\omega}^{-1}((N_d^{-1},1])\mid \omega\in \Omega\}$ 
is a cover of $V$. We may assume it to be generating. Hence
\[
Q^{p,q}_*(T,g,\WW,m)\le N_d^{m}\cdot \rho^{p,q}_*(T,g,\Phi,m)\quad\mbox{for $m\ge 1$}
\]
and, by Lemma~\ref{lm:Ki12},
\[
Q^{p,q}_*(T,g)=Q^{p,q}_*(T,g,\WW)\le N_d\cdot  \rho_*^{p,q}(T,g) 
=N_d\cdot  \rho^{p,q}(T,g) \, .
\]
We may apply this estimate to $T^k$ and  $g^{(k)}$ for $k\ge 1$.
Finally, we use both claims of (\ref{mul}) for large $k$ to obtain 
$
Q^{p,q}_*(T,g)\le \rho^{p,q}_*(T,g) 
$.
\end{proof}


\section{Spaces of anisotropic distributions and transfer operators on $\real^d$}\label{sec:spaces}
\label{S4}
In this section, we introduce Banach spaces of anisotropic distributions on $\real^d$ and then argue about the action of transfer operators on it. The argument in this section will be applied to iterates of our original diffeomorphism $T$ and weight $g$, using suitable local charts and partition of unity. 

For a subset $K\subset \real^d$ and $0\le s\le \infty$, 
let $C^s(K)$ be the set of $C^{s}$ functions $u:\real^d\to \complex$ whose support is contained in $K$.  Let  $C_0^{s}(\real^d)$ be the set of $u\in C^s(\real^d)$ with compact support. Let $C^{\infty}_*(\real^d)$ be the set of functions in $C^{\infty}(\real^d)$ such that $
\sup_{x \in \real^d} |\partial^\beta u(x)|<\infty
$
for all $\beta \in(\integer_+)^d$.
The Schwartz space $\SS$ consists of all 
 $u\in C^\infty(\real^d)$ that are {\em rapidly decaying}, that is, 
$
\sup_{x \in \real^d} |x^\alpha \partial^\beta u(x)|<\infty
$
for all $\alpha$, $\beta \in(\integer_+)^d$. So $C^{\infty}_0(\real^d)\subset \SS\subset 
C^{\infty}_*(\real^d)\subset C^{\infty}(\real^d)$.
\subsection{Definition of local spaces}\label{ss:localspaces}
The basic idea in the definition of our anisotropic spaces of distribution is to  slightly 
modify the classical Littlewood-Paley dyadic decomposition of functions in Fourier space, by 
introducing some cones of directions, or
"polarizations". This approach was introduced in our previous paper \cite{BT}. 
Below we modify the definitions in \cite{BT} slightly in order to get
the improved bounds in Theorem~\ref{mainprop}. (See also \cite{FR} for a recent Fourier analysis approach in
the analytic setting.)

For two cones $\cone$ and $\cone'$
in $\real^d$, we write $\cone \cc \cone'$  if
$\overline \cone\subset \mbox{ interior} \, (\cone' )\cup \{0\}$.
Let  $\cone_+$ and $\cone_-$ be
closed cones in $\real^d$  with nonempty interiors.  
Assume that  $\cone_+\cap\cone_-=\{0\}$ and  that $\cone_+$ and $\cone_-$ 
contain some $d_s$- and $d_u$-dimensional subspaces, respectively. 
Let 
$\varphi_+, \varphi_-:\sphere\to [0,1]$ be
$C^{\infty}$ functions on the unit sphere $\sphere$ in $\real^d$ satisfying
\[
\varphi_+(\xi)=\begin{cases}
1, &\mbox{if $\xi\in \sphere\cap \cone_{+}$,}\\
0, &\mbox{if $\xi\in \sphere\cap \cone_{-}$,}
\end{cases} \qquad 
\varphi_-(\xi)=1-\varphi_+(\xi).
\]
We shall work with combinations
$\Theta=(\cone_+,\cone_-,\varphi_+,\varphi_-)$ as above, which are called {\em polarizations.}

To a polarization $\Theta$ as above, we associate   the set $\FF=\FF(\Theta)$ 
of all $C^1$-submanifolds $F\subset \real^d$, of dimension $d_u$,
so that the straight line connecting any two distinct points in 
$F$ is normal to a $d_s$-dimensional subspace contained in  $\cone_+$.

\begin{remark}
Our assumption on $\FF(\Theta)$ 
implies that, if we take a $d_s$-dimensional subspace $E$  that is normal to 
a $d_u$-dimensional subspace $E^{\perp}\subset \cone_-$, then the projection $\pi: \real^d\to \real^d/E$ is a $C^1$ 
diffeomorphism when restricted to $F\in \FF(\Theta)$. 
\end{remark}

For $u\in C^\infty(\real^d)$ and $\FF=\FF(\Theta)$, we set
\begin{equation}
\|u \|_{L^1(\FF)}=\sup_{F\in \FF}
\|u\|_{L^1(\mu_F)}\in \real\cup \{\infty\} \, ,
\end{equation}
where $\mu_F$ is the Riemann volume on $F$ induced by the standard metric on $\real^d$.

The following lemma will play the role that the usual Young inequality for convolution played 
in \cite{BT}:
\begin{lemma}\label{psYoung}
Let $\FF=\FF(\Theta)$.
Then we have
\begin{align*}
\|A*u\|_{L^1(\FF)}&\le \|A\|_{L^1}\|u\|_{L^1(\FF)} \quad\mbox
{for $u\in C^\infty(\real^d)$ and $A\in L^1(\real^d)$}
\end{align*}
where $*$ denotes the convolution, 
$A*u(x)=\int_{\real^d} A(y) u (x-y) dy$.
\end{lemma}
\begin{proof}
Take $F\in \FF$ arbitrarily and let $F+x$ be the translation of $F$ by $x\in \real^d$, 
which also belongs to $\FF$. 
Then we have
\begin{align*}
\|(A*u)\|_{L^1(\mu_F)}&
\le \int_F\left( \int_{\real^d} |A(y)| \cdot |u(x-y)| dy\right) \, d\mu_F(x)\\
&= 
\int_{\real^d}|A(y)| \left( \int_{F}  |u(x-y)| d\mu_F(x)\right)\,  dy\\
&\le \int |A(y) |\cdot \|u\|_{L^1(\mu_{F-y})} \, dy
\le \|A\|_{L^1}\|u\|_{L^1(\FF)}\, ,
\end{align*}
where we used that $\mu_{F-y}$ is a translation of $\mu_F$.
\end{proof}

\medskip

We next introduce some notation in view of performing a
dyadic decomposition in the Fourier space. Let $\Theta=(\cone_+,\cone_-,\varphi_+,\varphi_-)$ be a polarization. 
Fix  a $C^\infty$ function $\chi:\real\to [0,1]$ with
$\chi(s)=1$ for $s\le 1$, and $\chi(s)=0$ for $s\ge 2$.
For $n\in \mathbb Z_+$, define 
$\chi_n:\real^d\to [0,1]$ by $\chi_n(\xi)=\chi(2^{-n}|\xi|)$,
and put $\chi_{-1}\equiv 0$. Set
$\psi_n:\real^d\to [0,1]$ to be 
$\psi_n(\xi)=\chi_n(\xi)-\chi_{n-1}(\xi)$, for $n \in \integer_+$.
Let $\Gamma=\{ (n,\sigma)\mid n\in \integer_+ ,\sigma\in\{+,-\}\}$.
For $(n,\sigma)\in \Gamma$,  we define
\[
\psi_{\Theta, n, \sigma}(\xi)=\begin{cases}
\psi_n(\xi)\varphi_{\sigma}(\xi/|\xi|),\quad&\mbox{ if $n\ge 1$,}\\
\psi_0(\xi)/2=\chi_0(\xi)/2,\quad&\mbox{ if $n=0$.}
\end{cases}
\]
Then the family of functions $\{\psi_{\Theta, n, \sigma}\}_{(n,\sigma)\in \Gamma}$ is a $C^{\infty}$ partition of unity. 
Note that the inverse Fourier transform
$
\widehat \psi_{\Theta, n,\sigma}(x)=(2\pi)^{-d}\int e^{ix\xi} \psi_{\Theta, n,\sigma}(\xi) \, d\xi
$
of each $\psi_{\Theta, n,\sigma}$ belongs to $\SS$, and satisfies the following scaling law:
\begin{equation}\label{eqn:scale}
\widehat \psi_{\Theta, n,\sigma}(x)=2^{d(n-1)} \widehat \psi_{\Theta, 1,\sigma}(2^{n-1} x) \quad \mbox{ if $n\ge 2$.}
\end{equation}
In particular, we have
\begin{equation}\label{eqn:invF}
\sup_{(n,\sigma)\in \Gamma} \|\widehat \psi_{\Theta, n,\sigma}\|_{L^1(\real^d)}<\infty.
\end{equation}

We may decompose $u\in C^{\infty}_0(\real^d)$ as
$u=\sum_{(n,\sigma)\in \Gamma}u_{\Theta,n,\sigma}$, by setting
\[
u_{\Theta,n,\sigma}:=\psi_{\Theta,n,\sigma}(D) u
=\widehat \psi_{\Theta,n,\sigma} * u \,\in \SS\, . 
\]
\begin{remark}
For $\psi\in \SS$, we define 
the pseudodifferential operator $\psi(D):\SS\to \SS$ by 
\[
\psi(D) u(x):=(2\pi)^{-d} \int e^{i\xi(x-y)} \psi(\xi) u(y) d\xi dy=\widehat{\psi}*u(x).
\]
We may write this operation as $\psi(D)=\Fourier^{-1}\circ M_{\psi}\circ \Fourier$ using Fourier transform $\Fourier$ and the multiplication operator $M_{\psi}$ by $\psi$. 
From the expression as a convolution operator, we may extend it as an operator 
$\psi(D):C^{\infty}_*(\real^d)\to C^{\infty}_*(\real^d)$.
We will often use the fact that
\begin{equation}\label{eqn:psed}
\psi_1(D)\psi_2(D)=(\psi_1\cdot \psi_2)(D),\quad 
\widehat{\psi}_1*\widehat{\psi_2}=\widehat{\psi_1\cdot \psi_2}
\quad \mbox{for $\psi_1,\psi_2\in \SS$.}
\end{equation}
\end{remark}

Here we quote the following lemma from \cite{BT}, which tells roughly that the functions $\psi_{\Theta,n,\sigma}(D)u$ decay rapidly outside of the support of $u\in C^{\infty}_0(\real^d)$.
\begin{lemma}[{\cite[{Lemma 4.1}]{BT}}]\label{lm:decay}
Let $K\subset \real^d$ be a compact subset. 
For any positive numbers $b,c$ and $\epsilon$, there exists a constant $C>0$ such that 
\begin{equation}\label{eqn:decay}
|\psi_{\Theta,n,\sigma}(D) u(x)|\le C\cdot \frac{\sum_{(\ell,\tau)\in \Gamma} 2^{-c\max\{n,\ell\}} \|\psi_{\Theta,\ell,\tau}(D)u\|_{L^{\infty}}}{d(x,\supp(u))^b}
\end{equation}
for any $(n,\sigma)\in \Gamma$, $u\in C^{\infty}(K)$ and $x\in \real^d$ satisfying $d(x,\supp(u))>\epsilon$. 
\end{lemma}
\begin{remark}\label{rem:decay}
Since 
$\psi_{\Theta,\ell,\tau}(D)u={\chi}_{\ell+1}(D) \psi_{\Theta,\ell,\tau}(D)u=\widehat{\chi}_{\ell+1}* \psi_{\Theta,\ell,\tau}(D)u$, we have 
\[
\|\psi_{\Theta,\ell,\tau}(D)u\|_{L^{\infty}}\le C  2^{d\ell} \|\psi_{\Theta,\ell,\tau}(D)u\|_{L^1(\FF(\Theta))}\; \mbox{ for any $(\ell,\tau)\in \Gamma$ and $u\in C^{\infty}_0(\real^d)$.}
\]
Therefore we may replace the  $L^\infty$ norm in (\ref{eqn:decay}) by the norm $\|\cdot\|_{L^1(\FF(\Theta))}$. 
\end{remark}

For a polarization $\Theta$, real numbers $q<0<p$ and $u\in C^{\infty}_0(\real^d)$, we define
\begin{equation}\label{def:normOnR}
\|u\|_{\BB^{\Theta,p,q}}=\max\biggl \{\; 
\sup_{n\ge 0}\; 2^{pn}\|  u_{\Theta,n,+}\|_{L^1(\FF(\Theta))}, 
\; \;\sup_{n\ge 0} \; 2^{qn}\|u_{\Theta,n,-}\|_{L^1(\FF(\Theta))}\;\biggr \}\, .
\end{equation}
Consider a non-empty compact subset $K\subset \real^d$. 
We first check that the definition above gives a norm on $C^{\infty}(K)$. 
Let $\|\cdot\|_{C^s}$ be the usual $C^s$ norm on $C^{s}(K)$. 
\begin{lemma}  \label{lm:Cs}
For any $s> p$, there exists a constant $C=C(s,K)$ such that $\|u\|_{\BB^{\Theta,p,q}} \le C \|u\|_{C^s}$ for all $u\in C^\infty(K)$. 
\end{lemma}
\begin{proof}
We may assume that $s$ is not an integer.
Recall the following characterization of $C^{s}$ norm in terms of Littlewood-Paley decomposition (see \cite[Appendix A]{Ta}): For non-integer $s>0$, the $C^s$ norm is equivalent to the norm defined by
\[
\|u\|_{C^s_*}:= \sup_{n\ge 0}\left(2^{sn} \|\psi_n(D)u\|_{L^\infty}\right)\,.
\]
Since $\psi_{\Theta,n,\sigma}(D)u=\sum_{m:|m-n|\le 1}\widehat{\psi}_{\Theta,m,\sigma}* (\psi_{n}(D) u)$ by (\ref{eqn:psed}),  we have 
\begin{equation}\label{eqn:pT}
\|\psi_{\Theta,n,\sigma}(D)u\|_{L^\infty}
\le C \| \psi_{n}(D) u\|_{L^\infty}
\quad\mbox{for any $(n,\sigma)\in \Gamma$}
\end{equation} 
by Young inequality and  (\ref{eqn:invF}).
Using Lemma \ref{lm:decay} with (\ref{eqn:pT}), we estimate $\psi_{\Theta,n,\sigma}(D)u$ outside some neighborhood of $K$ and obtain
\[
\|\psi_{\Theta,n,\sigma}(D)u\|_{L^1(\FF)}\le C(s,K)\cdot  2^{-s n} \|u\|_{C^s_*}\quad\mbox{for any $(n,\sigma)\in \Gamma$}\; .
\]
Clearly this implies the lemma.
\end{proof}

We may now give the definition of our anisotropic space  of distributions.
\begin{definition}
For a polarization $\Theta=(\cone_+,\cone_-,\varphi_+,\varphi_-)$
and  
real numbers $q<0<p$,  set $\BB^{\Theta,p,q}(K)$ to
be the completion of $C^{\infty}(K)$ (or, equivalently, that of $C^{s}(K)$ with $s>p$) with respect to  
$\|\cdot \|_{\BB^{\Theta, p,q}}$. 
\end{definition}
\begin{remark}
The only difference between the space $C^{\Theta,p,q}_*(K)$ in our previous paper \cite{BT}  and 
the space $\BB^{\Theta,p,q}(K)$ 
in the present work is that, in \cite{BT}, the 
norm $\|\cdot\|_{L^1(\FF)}$ in the definition above 
was the $L^\infty$ norm.  
\end{remark}

\begin{lemma}
\label{rem:sc}
For any $s>|q|$, the space $\BB^{\Theta,p,q}(K)$ is contained in the space of distributions of order $s$ supported on $K$. 
\end{lemma}
\begin{proof}
We may assume that $s$ is not an integer. 
Take any $u\in C^{\infty}(K)$ and $v\in C^{\infty}_0(\real^d)$ and decompose them as 
$u=\sum_{(n,\sigma)\in \Gamma} \psi_{\Theta,n,\sigma}(D)u$ and 
$v=\sum_{n\ge 0} \psi_n(D) v$ respectively. 
Since $\supp(\psi_{\Theta,n,\sigma})\cap \supp(\psi_{m})\neq \emptyset$ only if $|m-n|\le 1$, we get
\[
\int u\cdot v \,dx =\sum_{(n,\sigma)\in \Gamma} \sum_{m:|m-n|\le 1}\int \psi_{\Theta,n,\sigma}(D)u(x)  \cdot \psi_{m}(D)v(x)\, dx
\]
by Parseval's identity.
Using Lemma \ref{lm:decay} with Remark \ref{rem:decay}, we estimate $\psi_{\Theta,n,\sigma}(D)u$ outside some neighborhood of $K$ and obtain 
\[
\left|\int u\cdot v\, dx\right|\le C \|u\|_{\BB^{\Theta,p,q}}
\|v\|_{C^s_*}.
\]
This implies the claim of the lemma. 
\end{proof}

The decomposition introduced above can be viewed as an operator
\[
\QQ_{\Theta}:C^{\infty}(K)\to \SS^{\Gamma},\quad u\mapsto \bigl(u_{\Theta,n,\sigma}:=\psi_{\Theta,n,\sigma}(D)u\bigr)_{(n,\sigma)\in \Gamma}.
\]
Below we set up some Banach spaces for the target of $\QQ_{\Theta}$, in the place of $\SS^\Gamma$ above. For an integer $n\ge 0$, we define
\[
\BBf^{\Theta}_{n}=\{ u\in C^{\infty}(\real^d) \mid \chi_{n}(D)u=u 
\mbox{ and } \|u\|_{L^1(\FF(\Theta))}<\infty\}\,\subset C^{\infty}_*(\real^d)\, .
\]
For each $s\ge 0$ and $n\ge 0$, there exists a constant $C(s,n)>0$ such that 
\begin{equation}\label{eqn:cr}
\|u\|_{C^s}\le C(s,n)\|u\|_{L^1(\FF(\Theta))}
\qquad\mbox{for any $u\in \BBf^{\Theta}_{n}$,}
\end{equation}
because $\partial^\alpha u=(\partial^\alpha\widehat{\chi}_n)*u$.
Hence $\BBf^{\Theta}_{n}$ is a Banach space with respect to the norm 
$\|\cdot \|_{L^1(\FF(\Theta))}$. 
We have
that $\BBf^{\Theta}_{n-1}\subset \BBf^{\Theta}_{n}$ by $\chi_{n}(D)\chi_{n-1}(D)=\chi_{n-1}(D)$.

\begin{definition}
For a polarization
$\Theta$ and real
numbers $q<0<p$, we define
\[
\BBf^{\Theta,p,q}_{\Gamma}=\left\{(u_{n,\sigma})_{(n,\sigma)\in \Gamma}
\;\left|\; u_{n,\sigma}\in \BBf^{\Theta}_{n+3}, \;
\lim_{n\to \infty} \max_{\sigma=+,-}2^{c(\sigma)n}\|u_{n,\sigma}\|_{L^1(\FF(\Theta))}=0\right.\right\}\, ,
\]
where $c(+)=p$ and $c(-)=q$. This is a Banach space with respect to the norm
\begin{align*}
\|(u_{n,\sigma})_{(n,\sigma)\in \Gamma}\|_{\BBf^{\Theta, p,q}_\Gamma}
:=\sup_{(n,\sigma)\in \Gamma}
\left(
2^{c(\sigma)n} \|u_{n,\sigma}\|_{L^1(\FF(\Theta))}
\right)
 \, .
\end{align*}
\end{definition}
\begin{remark}\label{rem:B}
The space $\BBf^{\Theta,p,q}_{\Gamma}$ above is a closed subspace of the Banach space 
\[
\widetilde{\BBf}^{\Theta,p,q}_{\Gamma}=\left\{(u_{n,\sigma})_{(n,\sigma)\in \Gamma}
\;\left|\; u_{n,\sigma}\in \BBf^{\Theta}_{n+3}\mbox{ and }
\|(u_{n,\sigma})_{(n,\sigma)\in \Gamma}\|_{\BBf^{\Theta, p,q}_\Gamma}<\infty
\right.\right\}\, 
\]
with the identical norm $\|\cdot \|_{\BBf^{\Theta, p,q}_\Gamma}$.
The space $\BBf^{\Theta,p,q}_{\Gamma}$ is slightly more convenient than $\widetilde{\BBf}^{\Theta,p,q}_{\Gamma}$ for us. For instance the subset 
\begin{equation}\label{eqn:BZ}
\BBf_\Gamma^{\Theta}:=\{ (u_{n,\sigma})_{(n,\sigma)\in \Gamma} \mid u_{n,\sigma}\in \BBf^{\Theta}_{n+3}, \;\#\{(n,\sigma)\in \Gamma\mid u_{n,\sigma}\neq 0 \}<\infty \}
\end{equation}
is dense in $\BBf^{\Theta,p,q}_{\Gamma}$ though this is not true for $\widetilde{\BBf}^{\Theta,p,q}_{\Gamma}$. The difference will also make sense in Proposition \ref{nucltrace} and its proof.
\end{remark}

By (\ref{eqn:psed}), we have, for $k\ge 1$,  
\begin{equation}\label{eqn:chipsi}
\chi_{n+k}(D) \psi_{\Theta, n, \sigma}(D) 
=\psi_{\Theta, n, \sigma}(D) \quad \mbox{ on $C^{\infty}_*(\real^d)$.}
\end{equation}
This and Lemma \ref{lm:Cs} imply that $\QQ_{\Theta}(C^{\infty}(K))\subset \BBf^{\Theta,p,q}_{\Gamma}$. Thus, by the definitions of the norms, 
the operator $\QQ_{\Theta}$ above extends to the isometric embedding
\[
\QQ_{\Theta}: \BB^{\Theta, p,q}(K)\to  \BBf^{\Theta,p,q}_{\Gamma}\,  .
\]

From Lemma \ref{lm:decay}, we can see that the image of the embedding $\QQ_{\Theta}$ above is contained in much smaller subspaces than $\BBf^{\Theta,p,q}_{\Gamma}$. Indeed we can take a smaller Banch space $\hBBf^{\Theta,p,q}_{\Gamma} \subset \BBf^{\Theta,p,q}_{\Gamma}$ that contains the image of $\QQ_{\Theta}$ as follows. We set $
\beta(x)=(1+|x|^{2})^{(d+1)/2}$ and, for $n\ge 0$,
\[
\hBBf^{\Theta}_{n}=\{ u\in C^{\infty}(\real^d) \mid \chi_{n}(D)u=u 
\mbox{ and } \|\beta\cdot u\|_{L^1(\FF(\Theta))}<\infty\}\, .
\]
In parallel to (\ref{eqn:cr}), there exists a constant $C(s,n)>0$ for each $s\ge 0$ and $n\ge 0$ such that 
\begin{equation}\label{eqn:cr2}
\|\beta\cdot u\|_{C^s}\le C(s,n)\|\beta\cdot u\|_{L^1(\FF(\Theta))}
\qquad\mbox{ for any $u\in \BBf^{\Theta}_{n}$.}
\end{equation}
In particular, $\hBBf^{\Theta}_{n}$ is a Banach space with respect to the norm 
$u\mapsto \|\beta u \|_{L^1(\FF(\Theta))}$. 
\begin{definition}
For a polarization
$\Theta$ and two real
numbers $q<0<p$, we define
\[
\hBBf^{\Theta,p,q}_{\Gamma}=\left\{(u_{n,\sigma})_{(n,\sigma)\in \Gamma}
\;\left|\; u_{n,\sigma}\in \BBf^{\Theta}_{n+2},\;
\lim_{n\to \infty} \max_{\sigma=+,-}2^{c(\sigma)n}\|\beta\cdot u_{n,\sigma}\|_{L^1(\FF(\Theta))}=0\right.\right\}
\]
where $c(+)=p$ and $c(-)=q$. This is a Banach space with respect to the norm
\begin{align*}
\|(u_{n,\sigma})_{(n,\sigma)\in \Gamma}\|_{\hBBf^{\Theta, p,q}_\Gamma}
:=\sup_{(n,\sigma)\in \Gamma}
\left(
2^{c(\sigma)n} \|\beta\cdot u_{n,\sigma}\|_{L^1(\FF(\Theta))}
\right)
 \, .
\end{align*}
\end{definition}
Obviously the inclusion $\iota:\hBBf^{\Theta,p,q}_{\Gamma}\to \BBf^{\Theta,p,q}_{\Gamma}$ is non-expansive. 

\begin{lemma}\label{lm:hBBf} $\QQ_{\Theta}(\BB^{\Theta, p,q}(K))\subset \hBBf^{\Theta,p,q}_{\Gamma}$. $\QQ_{\Theta}:\BB^{\Theta, p,q}(K)\to \hBBf^{\Theta,p,q}_{\Gamma}$ is bounded.
\end{lemma}
\begin{proof}
By (\ref{eqn:chipsi}) for $k=2$ and by Lemma~\ref{lm:decay} with Remark \ref{rem:decay}, we can see that 
${\psi}_{\Theta, n,\sigma}(D) (C^{\infty}(K))\subset \BBf^{\Theta}_{n+2}$ and that $
\|Q_{\Theta}u\|_{\BBf^{\Theta,p,q}_{\Gamma}}\le \|Q_{\Theta}u\|_{\hBBf^{\Theta,p,q}_{\Gamma}}\le C\|Q_{\Theta} u\|_{\BBf^{\Theta,p,q}_{\Gamma}}$ 
for all $u\in C^{\infty}(K)$, for some constant $C$. This implies the lemma. 
\end{proof}


\subsection{Transfer operators associated to cone-hyperbolic maps}
\label{defs}

In this subsection, we define regular cone-hyperbolic maps on
bounded open subsets of $\real^d$
and consider transfer operators associated to such
maps $\TT$ and $C^{r-1}$ weights $G$. 

\begin{definition} Let $U$ and $U'$ be bounded open subsets in $\real^d$, and let 
$\Theta=(\cone_{+}, \cone_{-},\varphi_{+},\varphi_-)$ and $\Theta'=(\cone'_{+}, \cone'_{-},\varphi'_{+},\varphi'_-)$ be two 
polarizations\footnote{
We view $\cone_{\omega,\pm}$, $\cone'_{\omega,\pm}$ as constant cone fields in the {\em cotangent} bundle 
$T^*\real^d$, so we apply the {\em transpose} of $D\TT$ to the vectors in them.}.
A $C^r$ 
diffeomorphism  $\TT:U'\to U$ is {\em regular cone-hyperbolic} with respect to polarizations 
$\Theta$ and $\Theta'$  
if $\TT$ extends to a  bilipschitz $C^1$
diffeomorphism of $\real^d$ 
so that $D\TT_{x}^{tr}(\real^d\setminus \cone_{+}) \cc \cone'_{-}$ 
for each $x\in \real^d$
and, in addition, that there exists, for each $x,y \in \real^d$,
a linear transformation $\mathbb L_{xy}$ satisfying 
$(\mathbb L_{xy})^{tr}(\real^d\setminus \cone_{+}) \cc \cone'_{-}$
and 
$\mathbb L_{xy}(x-y)=\TT(x)- \TT(y)$.
(We denote the transposed matrix of $A$ by $A^{tr}$.)
\end{definition}
If $\TT$ is regular cone-hyperbolic, then the extension 
$\TT$ to $\real^d$ maps each element of $\FF(\Theta')$ to an element of
$\FF(\Theta)$, from both conditions in the definition. 
\begin{remark}\label{regular}
The second condition on the extension
of $\TT$ in the definition above  does not follow from the first condition. 
For example, consider  a hyperbolic horseshoe map $\TT$, and let $U$ be a small neighbourhood 
of the entire invariant
horseshoe.
\end{remark} 

In the rest of this section, we consider the transfer operator
\begin{equation}\label{def:L}
L:C^{r-1}(U)\to C^{r-1}(U'),\quad Lu=G\cdot (u\circ \TT)
\end{equation}
associated to a regular cone-hyperbolic $C^r$ diffeomorphism $\TT:U'\to U$ with respect to polarizations $\Theta$ and $\Theta'$ as above and a $C^{r-1}$ weight $G\in C^{r-1}(U')$.  

We begin with a simple estimate on the operator norm of $L$ with respect to the norms $\|\cdot \|_{L^1(\FF(\Theta))}$ and $\|\cdot \|_{L^1(\FF(\Theta'))}$. 
Define
\[
|\det (D\TT|_{\cone'_+})|(x):=\inf_{L} |\det (D\TT|_L)|(x)\, 
\quad \mbox{for $x\in U'$},
\]
where $\inf_{L}$ denotes the infimum over all $d_u$-dimensional subspaces $L\subset \real^d$ with normal subspace contained in $\cone'_+$,
and $\det (D\TT|_L)$ is defined as for (\ref{defdet}). 
Then we have, for any $u\in C^{r-1}(\real^d)$,   
\begin{equation}\label{eqn:chofv}
\|Lu \|_{L^1(\FF(\Theta'))}\le 
\|G\|_{L^\infty} \cdot \sup_{\supp(G)} 
\bigl( 
|\det D\TT|_{\cone'_+}|^{-1}\bigr)
\cdot \|u\|_{L^1(\FF(\Theta))}\, .
\end{equation}

Fix real numbers $q<0<p$ satisfying $p-q<r-1$ henceforth. Below we will introduce an auxiliary operator $M:\BBf^{\Theta,p,q}_{\Gamma}\to \BBf^{\Theta',p,q}_{\Gamma}$ and show that the following diagram of bounded operators commutes, with $L$ an extension of (\ref{def:L}):
\begin{equation}\label{CD:LM}
\begin{CD}
\BBf^{\Theta,p,q}_{\Gamma}@>{M}>> \BBf^{\Theta',p,q}_{\Gamma}\\
@A{\QQ_{\Theta}}AA  @A{\QQ_{\Theta'}}AA \\
\BB^{\Theta,p,q}(\overline{U})@>{L}>> \BB^{\Theta',p,q}(\overline{U'})
\end{CD}
\end{equation}
The operator $M$ is an infinite matrix of operators, each of which describes the transition between "frequency bands" induced by $L$. 

We recall some definitions from \cite{BT}. 
We associate, to $\TT$ and $G$, two integers 
\[
h_{\max}^{+}=h_{\max}^{+}(\TT,G)\quad\mbox{ and } \quad 
h_{\min}^{-}=h_{\min}^{-}(\TT,G)
\]
by
\begin{align*}
&h_{\max}^{+}=\biggl [\log_2 \biggl (\sup_{x\in \supp(G)} \sup_{
\stackrel{\scriptstyle{\|\xi\|=1}}{D\TT^{tr}(\xi)\notin \cone'_-}}
\| D\TT^{tr}_x(\xi)\|\biggr) \biggr ]+6\\
& h^{-}_{\min}=\biggl [\log_2 \biggl (\inf_{x \in \supp(G)} \inf_{
\stackrel{\scriptstyle{\|\xi\|=1}}{\xi\notin \cone_+}}
\| D\TT^{tr}_x(\xi)\|\biggr) \biggr ]-6 \, .
\end{align*}
\begin{remark}
We will consider the situation $h_{\max}^+\ll 0\ll h_{\min}^{-}$ in application.
\end{remark}

This definition implies that, for $x\in \supp(G)$ and $\xi\in \real^d$, 
\begin{equation}\label{eqn:minmax1}
\begin{aligned}
&\|D\TT^{tr}_{x}(\xi)\|< 2^{h_{\max}^{+}-5}\|\xi\|
\qquad\mbox{if $D\TT^{tr}_x(\xi)\notin \cone'_{-}$, and }\\
&2^{h_{\min}^{-}+5}\|\xi\|< \|D\TT^{tr}_x(\xi)\|
\qquad\mbox{if $\xi\notin \cone_{+}$.}
\end{aligned}
\end{equation}
We next introduce the relation $\hookrightarrow=\hookrightarrow_{\TT,G}$ on $\Gamma$ as follows: Write  $(\ell, \tau) \hookrightarrow (n,\sigma)$,
for $(\ell,\tau), (n,\sigma)\in \Gamma$,
if  either 
\begin{itemize}
\item $(\tau,\sigma)=(+,+)$ and $n\le \ell+h_{\max}^+$, or 
\item $(\tau,\sigma)=(-,-)$  and 
$\ell+h_{\min}^{-}\le n$, or
\item $(\tau,\sigma)=(+,-)$ and
($n\ge h_{\min}^-$ or $\ell\ge - h_{\max}^+$).
\end{itemize}
Otherwise we write $(\ell, \tau) \not\hookrightarrow (n,\sigma)$.

Take a closed cone $\widetilde \cone_{+} \cc \cone_{+}$
 such that
\begin{equation}\label{conehyp2}
D\TT_{x}^{tr}(\real^d\setminus \widetilde \cone_{+}) 
\cc \cone'_{-}\quad\mbox{for $x\in \supp(G)$}
\end{equation}
and another closed cone\footnote{Actually $\widetilde \cone_{-}$ will not play any roll in the following. One may set $\widetilde \cone_{-}=\emptyset$.}
 $\widetilde \cone_{-}\cc \cone_{-}$. Let  
$\tilde \varphi_{+}$, $\tilde\varphi_{-}:\sphere\to [0,1]$ be
$C^{\infty}$ functions satisfying
\[
\tilde \varphi_{+}(\xi)=\begin{cases}
1, &\mbox{if $\xi\notin \sphere\cap \cone_{-}$;}\\
0, &\mbox{if $\xi\in \sphere\cap \widetilde \cone_{-}$,}
\end{cases} \quad 
\tilde\varphi_{-}(\xi)=\begin{cases}
0, &\mbox{if $\xi\in \sphere\cap \widetilde\cone_{+}$;}\\
1, &\mbox{if $\xi\notin \sphere\cap \cone_{+}$.}
\end{cases}
\]
Put $\tilde\psi_\ell(\xi)=\chi(2^{-\ell-1}\|\xi\|)-\chi(2^{-\ell+2}\|\xi\|)$
for $\ell\ge 1$,  and define, for $(\ell,\tau)\in\Gamma$, 
\[
\tilde \psi_{\Theta, \ell, \tau}(\xi)=\begin{cases}
\tilde{\psi}_\ell(\xi)\tilde \varphi_{\tau}(\xi/\|\xi\|),& \mbox{ if $\ell\ge 1$,}\\
\chi(2^{-1}\|\xi\|),& \mbox{ if $\ell=0$.}
\end{cases}
\]
Then  
$\tilde \psi_{\Theta,\ell,\tau}(\xi)=1$ for all
$\xi \in \supp (\psi_{\Theta,\ell,\tau})$ and  
 (\ref{eqn:invF}) holds with $\psi_{\Theta,n,\sigma}$ replaced by $\tilde{\psi}_{\Theta,n,\sigma}$.
Further, by modifying the cone $\widetilde \cone_{+}$ if necessary, we may assume that  
\begin{equation}
\label{eqn:minmax2}
2^{h_{\min}^{-}+4}\|\xi\|<\|D\TT^{tr}_x(\xi)\|
\qquad\mbox{for any $x\in \supp(G)$ and any $\xi\notin \widetilde \cone_{+}$.}
\end{equation}

From (\ref{eqn:minmax1}--\ref{eqn:minmax2}), there exists a constant $C(\TT,G)>0$ such that, if $(\ell, \tau) \not\hookrightarrow (n,\sigma)$ 
and $\max\{n,\ell\}\ge C(\TT,G)$ for $(\ell,\tau),(n,\sigma)\in \Gamma$, then we have, for all $x \in \supp(G)$, 
\begin{align}\label{lowerbd}
&d(\supp(\psi_{\Theta',n,\sigma}), 
D\TT_x^{tr}(\supp(\tilde \psi_{\Theta,\ell,\tau})))
\ge 2^{\max\{n,\ell\}-C(\TT,G)} \,.
\end{align}

For each $(\ell,\tau),(n,\sigma)\in \Gamma$, we define the operator $S^{\ell,\tau}_{n,\sigma}:\BBf_{\ell+3}^{\Theta}\to \BBf^{\Theta'}_{n+1}$ 
by   
\[
S^{\ell,\tau}_{n,\sigma} u:=\psi_{\Theta',n,\sigma}(D) \circ 
L\circ  \tilde{\psi}_{\Theta,\ell,\tau}(D)  u\, .
\]
We begin with defining the operator $M$ formally by  
\begin{equation}\label{eqn:defM}
M((u_{n,\sigma})_{(n,\sigma)\in \Gamma})=\left(\sum_{(\ell,\tau)\in \Gamma}
S^{\ell,\tau}_{n,\sigma}u_{\ell,\tau}
\right)_{(n,\sigma)\in \Gamma}\;.
\end{equation}
To check that this formal definition gives a bounded operator $M:\BBf^{\Theta,p,q}_{\Gamma}\to \BBf^{\Theta',p,q}_{\Gamma}$, we recall from \cite{BT} a few estimates on the operators $S^{\ell,\tau}_{n,\sigma}$. 
Define the positive-valued integrable function $b:\real^d\to \real$ by
\begin{equation}\label{convol}
b(x)=1\quad\mbox{ if $\|x\|\le 1$}, \qquad
b(x)=\|x\|^{-d-1}\quad\mbox{ if $\|x\|> 1$.}
\end{equation}
For $m>0$, we set 
\begin{equation}\label{eqn:bm}
b_m:\real^d\to \real,\qquad b_m(x)=2^{dm} \cdot b(2^m x)\,,
\end{equation}
so that  $\|b_m\|_{L^1}=\|b\|_{L^1}$. There exists a constant $C>0$ such that 
\begin{equation}\label{eqn:bconv}
b_n * b_m(x)\le C \cdot b_{\min\{n,m\}}(x)  \quad \mbox{ for any $x\in \real^d$ and any $n,m\ge 0$.}
\end{equation}
By (\ref{eqn:scale}), there exists a constant $C>0$ such that, for any $x\in \real^d$ and $(n,\sigma)\in \Gamma$,  
\begin{equation}\label{eqn:kerpsi}
|\widehat{\psi}_{\Theta, n,\sigma}(x)|<C \cdot b_{n}(x)\quad \mbox{and} \quad  
|\widehat{\tilde{\psi}}_{\Theta, n,\sigma}(x)|<C \cdot b_{n}(x)\,.
\end{equation}
\begin{lemma}[{\cite[(27)]{BT}}]\label{lm:bounded}
There exists a constant 
$C(\TT,G)\ge 1$, which may depend on $\TT$ and $G$, so that, if 
$(\ell, \tau) \not\hookrightarrow (n,\sigma)$ for  $(\ell, \tau), (n,\sigma)\in \Gamma$, then we have
\begin{align} \label{inte}
&|S^{\ell,\tau}_{n,\sigma} u(x)|
\le C(\TT,G)\cdot 2^{-(r-1)\max\{n,\ell\}}\int_{\real^d} 
 b_{\min\{n,\ell\}}(x-y) \cdot |u(\TT(y))| \, dy\,,
\end{align}
for any $u\in C^{\infty}_*(\real^d)$ and $x\in \real^d$.
\end{lemma}
The proof is just a few applications of 
integration by parts using the estimate (\ref{lowerbd}).
For convenience of the reader, we give the proof in the case when $r$ is an integer 
in Appendix~\ref{apd:r}.

\begin{lemma}\label{lm:bdd} 
There is a constant $C> 1$, which does not depend on $\TT$ nor $G$, so that,
for any $(\ell,\tau), (n,\sigma)\in \Gamma$ and  any $u\in \BBf_{\ell+3}^{\Theta}$, it holds
\begin{equation}\label{eqn:E1}
\|\beta\cdot S^{\ell,\tau}_{n,\sigma} u\|_{L^1(\FF(\Theta'))}\le
C \cdot \|G\|_{L^{\infty}} 
\cdot \sup_{\supp(G)} 
\bigl(|\det D\TT|_{\cone'_+}|^{-1}\bigr) \cdot \|u\|_{L^1(\FF(\Theta))}
 \, .
\end{equation}
Further there is a constant
$C(\TT, G)> 1$ so that,
if 
$(\ell, \tau) \not\hookrightarrow (n,\sigma)$ in addition, then
\begin{equation}\label{eqn:E2}
\|\beta\cdot S^{\ell,\tau}_{n,\sigma} u\|_{L^1(\FF(\Theta'))}\le 
C (\TT,G)\cdot 2^{-(r-1)\max\{n,\ell\}} \|u\|_{L^1(\FF(\Theta))}\, .
\end{equation}
\end{lemma}
\begin{proof}
Lemma~\ref{psYoung} and (\ref{eqn:invF}) give the estimate 
\[
\|\psi_{\Theta',\ell,\tau}(D) u\|_{L^1(\FF(\Theta'))}=
\|\widehat{\psi}_{\Theta',\ell,\tau}* u\|_{L^1(\FF(\Theta'))}
\le C
\|u\|_{L^1(\FF(\Theta'))}
\]
and the parallel estimate for $\tilde{\psi}_{\Theta,\ell,\tau}(D)$. 
The claim (\ref{eqn:E1}) with $\beta$ replaced by $1$ 
follows from these estimates and (\ref{eqn:chofv}). 
The claim (\ref{eqn:E2}) with $\beta$ replaced by $1$ follows from Lemma~\ref{psYoung}, 
\ref{lm:bounded} and (\ref{eqn:chofv}). To put back the factor $\beta$, use\footnote{For (\ref{eqn:E2}), use also the fact that there exists constants $C(\TT,G)<C'(\TT,G)$ such that the relation $(\ell,\tau)\hookrightarrow (n,\sigma)$ holds if $\sigma n-\tau \ell<C(\TT,G)$ and only if $\sigma n-\tau \ell<C'(\TT,G)$.} the following consequence of Lemma \ref{lm:decay} and Remark \ref{rem:decay}: For any $\epsilon,c,b>0$, there exists a constant $C>0$ such that 
\[
|S^{\ell,\tau}_{n,\sigma} u(x)|\le C\cdot d(x,\overline{U'})^{-b}\cdot \sum_{(n',\sigma')\in \Gamma} 2^{- c \max\{n,n'\}} \|S^{\ell,\tau}_{n',\sigma'} u\|_{L^1(\FF(\Theta')) }
\]
for all $x\in \real^d$ with $d(x,\overline{U'})>\epsilon$.
\end{proof}

\begin{remark}\label{rem:C}
The constant $C$ in (\ref{eqn:E1})  depends only on the polarization $\Theta'$ and the family of functions  $\{\widetilde{\psi}_{\Theta,\ell,\tau}\}_{(\ell,\tau)\in \Gamma}$. 
On the contrary, the constant $C(\TT,G)$ in  (\ref{eqn:E2}) and (\ref{inte}) depends heavily on $\TT$ and $G$.  
\end{remark} 

\begin{corollary}\label{cor:ext}
The formal definition (\ref{eqn:defM}) gives a bounded linear operator $M:\BBf^{\Theta,p,q}_{\Gamma}\to \hBBf^{\Theta',p,q}_{\Gamma}$. The  transfer operator $L$ extends boundedly to   $L:\BB^{\Theta,p,q}(\overline{U})\to
\BB^{\Theta',p,q}(\overline{U'})$.  The diagram (\ref{CD:LM}) commutes and $M(\BBf^{p,q}_\Gamma)\subset \QQ_{\Theta'}(\BB^{\Theta',p,q}(\overline{U'}))$.  
\end{corollary}
\begin{proof}
The first claim is an immediate consequence of Lemma \ref{lm:bdd} and the definition of the relation $\hookrightarrow$.
It is then easy to check that $M\circ \QQ_{\Theta}=\QQ_{\Theta'}\circ L$ on $C^{\infty}(K)$. 
Recalling that $\QQ_{\Theta'}: \BB^{\Theta',p,q}(\overline{U'})\to \BBf^{\Theta',p,q}_\Gamma$ is an isometric embedding, we get the second claim and the commutative diagram (\ref{CD:LM}). 
Since $M(\BBf^{\Theta}_\Gamma)\subset \QQ_{\Theta'}(C^{\infty}(\overline{U'}))$ for
$\BBf^{\Theta}_{\Gamma}$ defined in (\ref{eqn:BZ}), we get the last claim by density. 
\end{proof}

In view of the argument above, it is natural to  
 decompose the operator $M$ into
\[
M_b((u_{n,\sigma})_{(n,\sigma)\in \Gamma})=\left(\sum_{(\ell,\tau):(\ell,\tau)\hookrightarrow (n,\sigma)}
S^{\ell,\tau}_{n,\sigma}u_{\ell,\tau}
\right)_{(n,\sigma)\in \Gamma}
\]
and 
\[
M_c((u_{n,\sigma})_{(n,\sigma)\in \Gamma})=\left(\sum_{(\ell,\tau):(\ell,\tau)\not\hookrightarrow (n,\sigma)}
S^{\ell,\tau}_{n,\sigma}u_{\ell,\tau}
\right)_{(n,\sigma)\in \Gamma}\;.
\]
By the same argument as in the case of $M$, we can check that the above definitions in fact gives bounded operators $M_b, M_c:\BBf^{\Theta,p,q}_{\Gamma}\to \hBBf^{\Theta',p,q}_{\Gamma}$ and 
$M=M_b+M_c$. Using Lemma \ref{lm:bdd} and the definition of the relation $\hookrightarrow$ more carefully, we get\footnote{We refer the proof of \cite[Theorem 6.1]{BT} for details, though it will not be necessary.}
\begin{lemma}\label{lm:Mb} 
There exists a constant $C>0$, which does not depend on $\TT$ nor $G$ such that the operator norm of $M_b: \BBf^{\Theta,p,q}_{\Gamma}\to \hBBf^{\Theta',p,q}_{\Gamma}$ is bounded by 
\[
C\cdot \|G\|_{L^{\infty}}\cdot 
\biggl(\sup_{\supp(G)} |\det D\TT|_{\cone'_+}|^{-1}\biggr)\cdot 2^{\max\{p\cdot h_{\max}^+,q\cdot h_{\min}^-\}}.
\] 
\end{lemma}

\subsection{Approximation numbers}
\label{appr}
We shall study approximation numbers of the operator $M_c$ and show, in particular, that $M_c$ is compact. First we recall some basic definitions and facts about the approximation number from \cite{Pie}.
Suppose that $\cB$ and $\hB$ are Banach spaces. 
For $k\in \integer_+$, we define
the $k$-th approximation number of a bounded linear operator $\PP:\cB\to \hB$ by
\[
a_k(\PP) = \inf \{
\|\PP-F \|_{\mathcal B} : F: \mathcal B \to  \hB\, , \hbox{ rank}(F) < k\} \, .
\]
For $1\le t\le \infty$, let $\LL^{(a)}_t(\cB,\hB)$ be the set of bounded linear operators $\PP:\cB\to \hB$ such that $(a_k(\PP))_{k\in \integer_+}\in \ell^t(\integer_+)$. For each $\PP\in \LL^{(a)}_t(\cB,\hB)$, we set 
$\|\PP\|^{(a)}_t:=\|(a_k(\PP))\|_{\ell^t}$. 

Suppose that $\PP_1:\cB_1\to \cB_2$ and  $\PP_2:\cB_2\to \cB_3$ are bounded linear operators on Banach spaces. If $\PP_1\in \LL^{(a)}_t(\cB_1,\cB_2)$ (resp.~$\PP_2\in \LL^{(a)}_t(\cB_2,\cB_3)$), 
then $\PP_2\PP_1\in \LL^{(a)}_t(\cB_1,\cB_3)$ and we have 
\begin{equation}\label{prod2}
\|\PP_2 \PP_1 \|_{t}^{(a)} \le  \|\PP_2\|_{L(\cB_2,\cB_3)} \|\PP_1\|^{(a)}_t\;\;\; (\mbox{resp.}\;
\|\PP_2 \PP_1 \|_{t}^{(a)} \le  \|\PP_2\|^{(a)}_t \|\PP_1\|_{L(\cB_1,\cB_2)} \;)
\end{equation}
where $\|\PP\|_{L(\cB,\hB)}$ denotes the operator norm of a linear operator $\PP:\cB\to \hB$. 
For $t,t',s\in [1,  \infty]$ such that $1/t+1/t'=1/s$, 
there is a constant  $C>0$ so that
\begin{equation}\label{mult}
\|\PP_2  \PP_1\|_{s}^{(a)}\le
 C  \|\PP_2\|^{(a)}_{t'}\cdot
 \|\PP_1\|^{(a)}_{t}\quad\mbox{ for  $\PP_1\in \LL^{(a)}_t(\cB_1,\cB_2),\;\PP_2\in \LL^{(a)}_{t'}(\cB_2,\cB_3)$} \, .\end{equation}
The next lemma tells that the operators in $\LL^{(a)}_1(\cB,\hB)$ have nuclear representations:
\begin{lemma}[{\cite[Proposition 2.3.11]{Pie}}]\label{lm:Pie} 
There is a constant $C>0$ such that, if $\PP\in \LL^{(a)}_1(\cB,\hB)$, there exist sequences $v_i\in \hB$ and $v_i^*\in \mathcal{B}^*$, $i\in \integer_+$, such that $\PP=\sum_i v_i\otimes v_i^*$ and  
$\sum_i \|v_i\|_{\hB} \|v_i^*\|_{\cB}\le C\|\PP\|_{1}^{(a)}$.
\end{lemma}
\begin{remark}\label{re:app} 
We refer \cite[2.2--2.3]{Pie} for more explanation about approximation number. In particular, we refer 2.3.3 and 2.2.9 of \cite{Pie} for (\ref{prod2}) and  (\ref{mult}) respectively. 
\end{remark}

We now return to the operator $M_c$. 
\begin{proposition}\label{powerdecay}
The operator $M_c:
\BBf^{\Theta,p,q}_\Gamma\to \hBBf^{\Theta',p,q}_\Gamma$ in the last subsection
belongs to  $\LL^{(a)}_t(\BBf^{\Theta,p,q}_\Gamma, \hBBf^{\Theta',p,q}_\Gamma)$ for
any $t>d/(r-p+q-1)$, and is hence compact. 
\end{proposition}

For the proof, we prepare
the following  approximation lemma:
\begin{lemma}\label{finiterank}
Let $K\subset \real^d$ be a compact subset and let $\Theta$ be a polarization.
Let $s >0$ and $\epsilon>0$ be positive real numbers. Then there exists a constant
$C>0$ such that, for each $N>0$ and $(n,\sigma)\in \Gamma$ with $n<N$, there exists an
operator $
F_{n,\sigma,N}:C^{r-1}(K)\to \hBBf^{\Theta}_{n+2}$
 of rank at most 
$2^{(1+\epsilon) d N}$, so that, for any $u\in C^{r-1}(K)$,
\begin{equation}\label{psiF}
\|\beta \cdot (\psi_{\Theta,n,\sigma}(D) u- F_{n,\sigma,N}\, u)\|_{L^1(\FF(\Theta))}
\le  C2^{-s N} \|u\|_{L^1(\FF(\Theta))}\,.
\end{equation}
\end{lemma}

\begin{proof}
We may assume that $K\subset (-1,1)^d$.  
Let $\phi:\real^d\to [0,1]$ be a $C^\infty$ function 
 so that $\phi\equiv 1$ 
on $[-1,1]^d$ and $\phi\equiv 0$ on outside of $[-2,2]^d$, and put $\phi_{a}(\xi)=\phi(2^{-a} \xi)$ for $a>0$. Take arbitrary $\epsilon>0$ and consider  $N>0$ and $(n,\sigma)\in \Gamma$ with $n<N$.

For $u\in C^{r-1}(K)$, put 
$H(u)=\phi_{\epsilon N}\cdot
\psi_{\Theta,n,\sigma}(D) u=\phi_{\epsilon N}\cdot\widehat{ \psi}_{\Theta,n,\sigma}* u$. Since 
the distance between $K$ and $\supp(1-\phi_{\epsilon N})$ is greater than $2^{\epsilon N}-1$, 
there exists a constant $C_s>0$ for any $s>0$ so that  
\begin{equation}\label{eqn:upsi1}
\|\beta\cdot(H(u)- \psi_{\Theta,n,\sigma}(D)
u)\|_{L^1(\FF(\Theta))}\le C_s\cdot  2^{-sN}\|u\|_{L^1(\FF(\Theta))}\,\quad \mbox{for $u\in C^{r-1}(K)$.}
\end{equation}

Since $H(u)$ for $u\in C^{r-1}(K)$ is supported on $(-2^{\epsilon N+1}, 2^{\epsilon N+1})^d$, 
we may regard it as a function on $\real^d/(2^{\epsilon
N+2}\mathbb{Z})^d$ and consider the discrete Fourier coefficients
\begin{equation*}
c_\alpha (u)=(2^{\epsilon N+2})^{-d/2}
\int_{\real^d}
e^{-i \alpha x} H(u)(x) dx\, \qquad\mbox{for $\alpha \in (2^{-\epsilon N-1}\pi) \cdot  \integer^d$ . }
\end{equation*}
Set
\[
F(u)(x)=\phi_{\epsilon N+1}(x)\cdot
\sum_{|\alpha|\le 2^{N+5} }c_\alpha(u) \cdot e^{i\alpha x} \,\qquad \mbox{for  $u\in C^{r-1}(K)$} .
\]
Then the difference $H(u)- F(u)$ is supported on $(-2^{\epsilon N+2}, 2^{\epsilon N+2})^d$ and satisfies 
\begin{equation}\label{pdiff}
\|\beta\cdot (H(u) - F(u))\|_{L^\infty}\le
\left(\sup_{(-2^{\epsilon N+2}, 2^{\epsilon N+2})^d}\beta\right)\cdot 
\sum_{|\alpha| > 2^{N+5} }|c_\alpha(u)| \, .
\end{equation}
We may write the coefficient $c_{\alpha}(u)$ for $\alpha \in (2^{-\epsilon N-1}\pi) \cdot  \integer^d$ as
\[
(2^{\epsilon N+2})^{-d/2}\cdot \Fourier (H(u))(\alpha)
(2^{\epsilon N+2})^{-d/2}\cdot 
(\Fourier(\phi_{\epsilon N})* (\psi_{\Theta,n,\sigma}\cdot \Fourier(u)) )(\alpha)
\, 
\]
where $\Fourier$ denotes Fourier transfrom. 
We have that $\Fourier(\phi_{\epsilon N})(\xi)=2^{\epsilon N d} \cdot  \Fourier\phi(2^{\epsilon N}\xi)$ with $\Fourier\phi\in \SS$ and that $
\|\Fourier(u)\|_{L^\infty}
\le \|u\|_{L^1}
\le C(K) \|u\|_{L^1(\FF(\Theta))}$. 
Also we have  $|\alpha-\xi|>2^{-1}|\alpha|>2^{N+4}$ if
 $\xi \in \supp ( \psi_{\Theta,n,\sigma})$ and $|\alpha| > 2^{N+5}$.
Therefore, for any $s>0$, there exists a constant $C_s>0$ such that  
\[
|c_\alpha | \le C_s \cdot |\alpha|^{-s }\|u\|_{L^1(\FF(\Theta))}\quad\mbox{ if $|\alpha| > 2^{N+5}$.}
\]
Using this estimate in (\ref{pdiff}) and recalling (\ref{eqn:upsi1}), we find a constant $C_s>0$ for each $s>0$ so that  
\begin{equation}\label{eqn:pf}
\|\beta\cdot (\psi_{\Theta,n,\sigma}(D) u - F(u))\|_{L^1(\FF(\Theta))}\le  C_s \cdot 2^{-sN}\|u\|_{L^1(\FF(\Theta))}
\quad \mbox{ for  $u\in C^{r-1}(K)$.}
\end{equation}

Finally put $
F_{n,\sigma,N}(u):=\chi_{n+1}(D)(F(u))$ for $u\in C^{r-1}(K)$. The rank of the operator $F_{n,\sigma,N}$ or that of $F$ is bounded by
\[
\#\{ \alpha \in 2^{-\epsilon N-1}\pi \cdot  \integer^d\mid |\alpha|\le 2^{N+5}\}<C 2^{(1+\epsilon)dN}.
\]
It is not difficult to see that there exists a constant $C>0$ such that 
\[
\|\beta\cdot \chi_{n+1}(D)v\|_{L^1(\FF(\Theta))}
\le C\|\beta\cdot v\|_{L^1(\FF(\Theta))}\quad \mbox{for any $n\ge 0$ and $v\in C^\infty_0(\real^d)$.}
\]
Thus the claim (\ref{psiF}) follows from (\ref{eqn:chipsi}) and (\ref{eqn:pf}): 
\begin{align*}
\|\beta\cdot (\psi_{\Theta,n,\sigma}(D) u - F_{n,\sigma,N}(u))\|_{L^1(\FF(\Theta))}&\le  C\cdot 
\|\beta\cdot (\psi_{\Theta,n,\sigma}(D) u - F(u))\|_{L^1(\FF(\Theta))}
\\
&\le  C_s \cdot 2^{-sN}\|u\|_{L^1(\FF(\Theta))}\,.
\end{align*}
From (\ref{psiF}) and the relation $\chi_{n+2}(D)\chi_{n+1}(D)=\chi_{n+1}(D)$, the image of $F_{n,\sigma,N}$ is contained in $\hBBf^{\Theta}_{n+2}$. 
\end{proof}

\begin{proof}[Proof of Proposition \ref{powerdecay}]
We first approximate the operators $S^{\ell,\tau}_{n,\sigma}$ defined in the last subsection by finite rank operators.
By Lemma \ref{psYoung} and  (\ref{eqn:chofv}), we have
\[
\|L\circ \tilde{\psi}_{\Theta,\ell,\tau}(D) u\|_{L^1(\FF(\Theta'))}\le C(\TT,G)\cdot \|u\|_{L^1(\FF(\Theta))}
\quad \mbox{for any $(\ell,\tau)\in \Gamma$ and $u\in \BBf^{\Theta}_{\ell+3}$.}
\]
Take arbitrary $\epsilon>0$ and let $N>0$. Applying 
Lemma ~\ref{finiterank} to approximate the post-composition of $\psi_{\Theta',n,\sigma}(D)$, 
we find an operator $F_{n,\sigma}^{\ell,\tau}:
\BBf^{\Theta}_{\ell+3}\to \hBBf^{\Theta'}_{n+2}$ of rank at most $2^{(1+\epsilon)dN}$ for each $(n,\sigma), (\ell,\tau)\in \Gamma$ with $n<N$, such that
\begin{equation}\label{above}
\|\beta\cdot (S_{n,\sigma}^{\ell,\tau}(u)-
F_{n,\sigma}^{\ell,\tau}(u))\|_{L^1(\FF(\Theta'))}
\le C(\TT,G)\cdot  2^{-(r-1)N}\|u\|_{L^{1}(\FF(\Theta))}.
\end{equation}

Define $\mathbb{P}_N:\BBf^{p,q}_Z\to \hBBf^{p,q}_Z$ by 
\[
\mathbb{P}_N((u_{\ell,\tau})_{(\ell,\tau)\in \Gamma})=\left(\sum_{(\ell,\tau)\in \Gamma} P_{n,\sigma}^{\ell,\tau}(u_{\ell,\tau})
\right)_{(n,\sigma)\in \Gamma}
\]
where $P_{n,\sigma}^{\ell,\tau}=F_{n,\sigma}^{\ell,\tau}
$ if $\max\{n,\ell\}<N$ and $(\ell,\tau)\not\hookrightarrow (n,\sigma)$, and 
$P_{n,\sigma}^{\ell,\tau}=0$ otherwise. The rank of $\mathbb P_N$ is bounded by $C \cdot N^2 \cdot 2^{(1+\epsilon) d N}$.
By (\ref{above}) and the claim (\ref{eqn:E2}) of Lemma \ref{lm:bdd}, we obtain
\[
\|M_{c}-\mathbb P_N\|_{L(\BBf^{\Theta,p,q}_\Gamma, \hBBf^{\Theta',p,q}_\Gamma)}\le 
C(\TT,G,\epsilon)\cdot  2^{-(r-p+q-1-\epsilon)N}\, .
\]
This implies that $a_{k}(M_{c})<C 2^{-(r-p+q-1-\epsilon)N}$ for
$k=[C\cdot N^2\cdot 2^{(1+\epsilon) d N}]+1$, so that
$\|(a_{k}(M_{c}))\|_{ \ell^{t}}$ is bounded  for any $t>d(1+\epsilon)/(r-p+q-1-\epsilon)$. Since $\epsilon>0$ is arbitrary, we get the proposition. 
\end{proof}

\section{The transfer operator $\LL$ and its extensions}
\label{sec:L}
In this section, we study the transfer operator $\LL=\LL_{T,g}$ for a $C^r$ diffeomorphism $T:X\to X$ and a weight $g\in C^{r-1}(V)$, within the setting in Section~\ref{sec:intro}. 
Using local charts and a partition of unity, we associate to $\LL$ a system $\KK$ of transfer operators on local charts and then introduce a key auxiliary operator $\MM$. 
Once we define the operators $\KK$ and $\MM$ and check their relations to $\LL$, the proof of Theorem \ref{mainprop} is an immediate consequence of the argument in the last section. 

\subsection{Local charts adapted to the hyperbolic structure}
\label{localcharts}
We first set up a finite system of $C^{\infty}$ local charts 
on $V$, and of  polarizations on each of the local charts, so that they are adapted to the hyperbolic 
structure of the dynamical system $T$.
Consider $C^{\infty}$ local charts $\{(V_\omega, \kappa_\omega)\}_{\omega\in \Omega}$, with 
open subsets
$V_\omega\subset X$ and maps
$\kappa_\omega : V_ \omega\to \real^d$ such that $V \subset \cup_\omega V_\omega$, and  
consider also a system of polarizations
on those local charts $
\{\Theta_\omega=(\cone_{\omega,+}, \cone_{\omega,-},\varphi_{\omega,+}\, , 
\varphi_{\omega,-})\}_{\omega\in \Omega}$. 
Since $T$ is hyperbolic 
on $\Lambda$, we may assume that the following conditions hold:

\renewcommand{\labelenumi}{(\alph{enumi})}
\begin{enumerate}
\item $\mathcal{V}=\{V_\omega\}_{\omega\in \Omega}$ is a generating cover of $V$ and there is no strict subcover.
\item $U_{\omega}=\kappa_\omega(V_\omega)$  
is a bounded open subset of $\real^d$ for each $\omega\in \Omega$.
\item If $x\in V_\omega\cap \Lambda$, the cone
$(D\kappa_\omega)^{*}_x(\cone_{\omega,+})$    
contains the normal subspace of $E^u(x)$, and the cone
$(D\kappa_\omega)^{*}_x(\cone_{\omega,-})$
contains the normal subspace of  $E^s(x)$.
\item 
If $V_{\omega'\omega}=T^{-1}(V_\omega)\cap V_{\omega'}\ne \emptyset$,  the map in charts 
\[
T_{\omega' \omega}=\kappa_\omega
\circ T\circ \kappa_{\omega'}^{-1}
:\kappa_{\omega'}(V_{\omega'\omega}) \to U_\omega\, 
\]
is a $C^r$ regular cone-hyperbolic
diffeomorphism with respect to the polarizations 
$\Theta_\omega$ and $\Theta_{\omega'}$. 
\end{enumerate}
Let $\Phi=\{\phi_\omega\}$ be a $C^{\infty}$ partition of  unity for $V$ subordinate to the cover $\{V_\omega\}_{\omega\in \Omega}$, 
that is, the support of each $\phi_\omega:X\to [0,1]$ is contained in $V_\omega$, and 
we have $\sum_{\omega\in \Omega} \phi_\omega(x)= 1$ for all $x\in V$. 
We will henceforth fix the local charts, the system of polarizations and the 
partition of unity as above. 
We may now define the space $\BB^{p,q}(T,V)$ of distributions:
\begin{definition}\label{defspace}The 
Banach space $\BB^{p,q}(T,V)$ is
the completion of $C^{\infty}(V)$ for the norm
\[
\|\varphi\|_{\BB^{p,q}(T,V)}=\max_{\omega\in \Omega} 
\|(\phi_\omega\cdot \varphi)\circ \kappa_\omega^{-1}\|_{\BB^{\Theta_\omega, p,q}}\, 
\] 
where the norms $\|\cdot \|_{\BB^{\Theta_\omega, p,q}}$ are those defined by (\ref{def:normOnR}). 
\end{definition}
By Lemma \ref{lm:Cs} and \ref{rem:sc}, the space $\BB^{p,q}(T,V)$ contains $C^s(V)$ for each $s> p$ and contained in the dual of $C^s(X)$ for each
$s>|q|$.

We decompose the iterates $\LL^m$ of $\LL$ as follows. Take
\footnote{We need to take the function $\tilde{g}$ in treating the case $\inf |g|=0$. Otherwise we may set $\tilde{g}=g$.}  
a positive-valued $C^{r-1}$ function $\tilde{g}:X\to \real$ such that $\tilde{g}(x)>|g(x)|$ for $x\in X$. For each $m\ge 1$, choose a subset $\Omega_m \subset \Omega^{m}$ so that 
\[
\VV_m:=\{V_{\vomega}:=\cap_{i=0}^{m-1}T^{-i}(V_{\omega_i})\mid \vomega=(\omega_0, \omega_1,\dots, \omega_{m-1})\in \Omega_m\}
\]
is a cover of $\bigcap_{k=0}^{m-1}T^{-k}(V)$ by non-empty open sets and that  (recall (\ref{defQ_0}))
\[
Q^{p,q}_*(T,\tilde{g},\VV,m)=\sum_{\vomega\in \Omega_m} \sup_{V_{\vomega}}\left( 
\frac{|\tilde{g}^{(m)}|\cdot  \lambda^{(p,q,m)}}{|\det(DT^m|_{E^u})|}\right)\, .
\]
Take a $C^{\infty}$ partition of unity 
$\Phi_m=\{ \phi_{\vomega}\in C^{\infty}(V_{\vomega}) \mid \vomega\in \Omega_m\}$ for $\bigcap_{k=0}^{m-1}T^{-k}(V)$
subordinate to $\VV_m$. 
Then we have $\LL^m=\sum_{\vomega\in \Omega_m}\LL^m_{\vomega}$ for the operators 
\[
\LL^m_{\vomega}:C^{r-1}(V)\to C^{r-1}(V),\qquad  
\LL^m_{\vomega}\varphi=\phi_{\vomega}\cdot g^{(m)}\cdot \varphi \circ T^{m}\;.
\]

\subsection{The system of transfer operators on local charts} 
\label{system}
We introduce the operator $\KK$ as follows. 
For each $\omega\in \Omega$, take a $C^\infty$ function $h_{\omega}\in C^{\infty}(U_\omega)$ so that $0\le h_{\omega}\le 1$ on $\real^d$ and that $h_{\omega}\equiv 1$ on $\kappa_\omega(\supp(\phi_\omega))$.
Set $C^{r-1}_{\Omega}=\bigoplus_{\omega\in \Omega} C^{r-1}(\overline{U_{\omega}})$. 
We define the operators $\Phi_*:C^{r-1}(V)\to C^{r-1}_{\Omega}$ and $H:C^{r-1}_{\Omega}\to C^{r-1}(V)$  by 
\[
\Phi_*(u)=((\phi_{\omega}\cdot u)\circ \kappa_{\omega}^{-1})_{\omega\in \Omega}\quad\mbox{and}\quad
 \quad H(( u_{\omega})_{\omega\in \Omega})
=\sum_{\omega\in \Omega} (h_{\omega}\cdot u_{\omega})\circ \kappa_\omega\;. 
\]
Obviously we have $H\circ \Phi_*=\id$. 
For each $m\ge 1$, we define  
\[
\KK^m=\Phi_*\circ \LL^m \circ H:C^{r-1}_{\Omega}\to C^{r-1}_{\Omega}\;.
\]
\begin{remark}\label{MTO}
The operator $\KK^m$ can be regarded as a regular MTO in the sense of Kitaev \cite{Ki}. 
\end{remark}

Then $\KK^m$ is the $m$-th iterate of $\KK:=\KK^1$ and the 
following diagram commutes:
\begin{equation}\label{CD:hp}
\begin{CD}
C^{r-1}(V)@>{\Phi_*}>> C^{r-1}_{\Omega}\\
@V{\LL^m}VV  @V{\KK^m}VV \\
C^{r-1}(V)@>{\Phi_*}>> C^{r-1}_{\Omega}\, .
\end{CD}
\end{equation}
Likewise, for $m\ge 1$ and $\vomega\in \Omega_m$, we define the operator 
$\KK_{\vomega}^m$ by replacing $\LL$ by $\LL_\vomega^m$ in the definition of $\KK^m$. Then we have $\KK^m=\sum_{\vomega\in \Omega_m}\KK^m_{\vomega}$ and the commutative diagram above holds with $\LL^m$ and $\KK^m$ replaced by $\LL_\vomega^m$ and $\KK_\vomega^m$ respectively.

We can describe the operator $\KK_{\vomega}^m$ as follows. Set 
\[
U_{\vomega,\omega'\omega}:=\kappa_{\omega'}(V_{\omega'}\cap V_{\vomega}\cap  T^{-m}(V_{\omega}))
\]
and define $T^m_{\omega'\omega}:U_{\vomega,\omega'\omega}\to U_{\omega}$ and $G_{\vomega,\omega'\omega}\in C^{r-1}(U_{\vomega,\omega'\omega})$ by
\[
T^m_{\vomega,\omega'\omega}=\kappa_{\omega}\circ T^m\circ \kappa_{\omega'}^{-1},\quad
G_{\vomega, \omega'\omega}=((\phi_{\omega'}\cdot \phi_{\vomega}\cdot g^{(m)})\circ \kappa_{\omega'}^{-1}) \cdot 
(h_{\omega}\circ T^m_{\vomega,\omega'\omega}).
\] 
For $\omega, \omega'\in \Omega$, we define
\[
(\KK^m_{\vomega})_{\omega\omega'}:C^{r-1}(\real^d)
\to C^{r-1}(U_{\vomega,\omega'\omega}),\qquad
(\KK^m_{\vomega})_{\omega\omega'} u=G_{\vomega, \omega'\omega}\cdot (u\circ  
T^m_{\vomega,\omega'\omega}).
\]
Then these operators are $\omega\omega'$-components of $\KK_{\vomega}^m$:
\[
\KK^m_{\vomega}((u_{\omega})_{\omega\in \Omega})=\left(
\sum_{\omega\in \Omega} (\KK^m_{\vomega})_{\omega\omega'}u_{\omega}
\right)_{\omega'\in \Omega}\;.
\]

We will apply the argument in the last section to $L=(\KK_{\vomega}^m)_{\omega\omega'}$, setting 
\begin{gather}\label{setting}
\TT=T^m_{\vomega,\omega'\omega}, \; G=G_{\vomega, \omega'\omega}, \; \Theta=\Theta_\omega, \; \Theta'=\Theta_{\omega'},\;\\ U'=U_{\vomega,\omega'\omega},
\; U=T^m_{\vomega,\omega'\omega}(U_{\vomega,\omega'\omega}).\notag
\end{gather}
For this purpose, we have to choose cones $\widetilde \cone_{\omega,+} \cc \cone_{\omega,+}$, $\widetilde \cone_{\omega,-} \cc \cone_{\omega,-}$  for each $\omega\in \Omega$, so that, for any $m\ge 1$ and $\vomega\in \Omega_m$,  if we set 
\begin{equation}\label{setting2}
\widetilde \cone_{+}=\widetilde \cone_{\omega,+}\quad \mbox{ and }\quad \widetilde \cone_{-}=\widetilde \cone_{\omega,-}
\end{equation}
in addition to (\ref{setting}), the conditions (\ref{conehyp2}) and (\ref{eqn:minmax2}) hold. 
Clearly this is possible if we take $\widetilde \cone_{\omega,+}$ sufficiently close to $\cone_{\omega,+}$.
We then choose $C^{\infty}$ functions $\tilde{\varphi}_{\omega,+}, \tilde{\varphi}_{\omega,-}:\sphere\to [0,1]$ and define $\tilde{\psi}_{\Theta_\omega,n,\sigma}\in C^{\infty}_0(\real^d)$ in the way parallel to that in the definitions of  
$\tilde{\varphi}_{+}$, $\tilde{\varphi}_{-}$ and $\tilde{\psi}_{\Theta,n,\sigma}$ in Subsection \ref{defs}.
When we refer the setting (\ref{setting}) in the following, we understand that it includes the additional setting (\ref{setting2}) and 
\begin{equation}\label{setting3}
\tilde{\varphi}_{+}=\tilde{\varphi}_{\omega,+},\quad 
\tilde{\varphi}_{-}=\tilde{\varphi}_{\omega,-}\quad 
\mbox{and}\quad 
\tilde{\psi}_{\Theta,n,\sigma}=\tilde{\psi}_{\Theta_\omega,n,\sigma}\quad \mbox{ for $(n,\sigma)\in \Gamma$.}
\end{equation}

Consider the Banach space 
\[
\BB^{p,q}_\Omega= \bigoplus_{\omega\in \Omega} \BB^{\Theta_\omega,p,q}\bigl(\overline{U_\omega}\bigr)
\]
with the norm $
\|(u_\omega)_{\omega\in \Omega}\|_{\BB^{p,q}_\Omega}
=\max_{\omega\in \Omega} \|u_\omega\|_{\BB^{\Theta_\omega, p,q}}$. 
By the definitions of the norms, the operator $\Phi_*$ extends to an isometric embedding 
$\Phi_*:\BB^{p,q}(V,T)\to \BB^{p,q}_\Omega$. 
Corollary \ref{cor:ext} applied to the setting (\ref{setting}) tells that the diagram (\ref{CD:hp}) extends to the following commutative diagram of bounded operators:
\[
\begin{CD}
\BB^{p,q}(V,T)@>{\Phi_*}>> \BB^{p,q}_\Omega\\
@V{\LL^m_{\vomega}}VV  @V{\KK^m_{\vomega}}VV \\
\BB^{p,q}(V,T)@>{\Phi_*}>> \BB^{p,q}_\Omega
\end{CD}
\]
Taking the sum with respect to $\vomega$, we get the same commutative diagram with $\LL^m_{\vomega}$ and $\KK^m_{\vomega}$ replaced by $\LL^m$ and $\KK^m$.

\subsection{The auxiliary operator $\MM$}
\label{sec:MM}
We next introduce the auxiliary operator $\MM$ as follows. 
Recall the Banach spaces $\BBf^{\Theta,p,q}_{\Gamma}$ in the last section and consider the Banach spaces 
\[
\BBf_Z^{p,q}:=\bigoplus_{\omega\in \Omega} \BBf^{\Theta_\omega,p,q}_\Gamma,
\qquad\hBBf_Z^{p,q}:=\bigoplus_{\omega\in \Omega} \hBBf^{\Theta_\omega,p,q}_\Gamma
.
\]
with the norms
\[
\|(\mathbf{u}_\omega)_{\omega\in \Omega}\|_{\BBf^{p,q}_Z}
:=\max_{\omega\in \Omega} \|\mathbf{u}_\omega\|_{\BBf^{\Theta_\omega,p,q}_\Gamma} \,, 
\quad
\|(\mathbf{u}_\omega)_{\omega\in \Omega}\|_{\hBBf^{p,q}_Z}
:=\max_{\omega\in \Omega} \|\mathbf{u}_\omega\|_{\hBBf^{\Theta_\omega,p,q}_\Gamma}\;.
\]
Let $\QQ:\BB^{p,q}_\Omega\to \BBf_Z^{p,q}$ be the isometric embedding defined by
\[
\QQ((u_\omega)_{\omega\in \Omega})=(\QQ_{\Theta_\omega}(u_\omega))_{\omega\in \Omega}\;.
\]
Applying the construction in Subsection \ref{defs} to $L=(\KK_{\vomega}^m)_{\omega\omega'}$ in the setting (\ref{setting}), we define the operator
\[
M=(\MM^m_{\vomega})_{\omega\omega'}:\BBf^{\Theta_\omega,p,q}_\Gamma\to \hBBf^{\Theta_{\omega'},p,q}_\Gamma
\]
for $\vomega\in \Omega_m$ and $\omega,\omega'\in \Omega$, so that the following diagram commutes:
\begin{equation}\label{CD:Mww}
\begin{CD}
\BB^{\Theta_\omega,p,q}\bigl(\overline{U_\omega}\bigr)@>{\QQ_{\Theta_\omega}}>> \BBf^{\Theta_\omega,p,q}_\Gamma\\
@V{(\KK^m_{\vomega})_{\omega\omega'}}VV  @V{(\MM^m_{\vomega})_{\omega\omega'}}VV \\
\BB^{\Theta_{\omega'},p,q}\bigl(\overline{U_{\omega'}}\bigr)@>{\QQ_{\Theta_{\omega'}}}>> \BBf^{\Theta_{\omega'},p,q}_\Gamma\end{CD}
\end{equation}
We define the bounded operator $\MM^m_{\vomega}:\BBf^{p,q}_Z\to \BBf^{p,q}_Z$ by
\[
\MM^m_{\vomega}((\mathbf{u}_\omega)_{\omega\in \Omega})
= \left(
\sum_{\omega\in \Omega} (\MM^m_{\vomega})_{\omega\omega'}\mathbf{u}_\omega
\right)_{\omega'\in \Omega}\;
\]
and put $\MM^m=\sum_{\vomega\in \Omega_m} \MM^m_{\vomega}$. Then we obtain the following commutative diagram of bounded operators:
\begin{equation}\label{CD:full}
\begin{CD}
\BB^{p,q}(T,V)@>{\Phi_*}>>\BB^{p,q}_\Omega@>{\QQ}>> \BBf^{p,q}_Z\\
@V{\LL^m_{\vomega}}VV  
@V{\KK^m_{\vomega}}VV  @V{\MM^m_{\vomega}}VV \\
\BB^{p,q}(T,V)@>{\Phi_*}>>\BB^{p,q}_\Omega@>{\QQ}>> \BBf^{p,q}_Z\end{CD}
\end{equation}
and the same diagram with $\LL^m_{\vomega}$, $\KK^m_{\vomega}$ and $\MM^m_{\vomega}$ replaced by 
$\LL^m$, $\KK^m$ and $\MM^m$. 

By using continuity of $\MM^m$, we can check that $\MM^m$ is the $m$-th iteration of $\MM:=\MM^1$ and that\footnote{See the proof of Corollary \ref{cor:ext} for the second inclusion.}
\[
\KK^m(\BB^{p,q}_\Omega)\subset \Phi_*(\BB^{p,q}(T,V)) \quad\mbox{and}\quad
\MM^m(\BBf^{p,q}_Z)\subset \QQ(\BB^{p,q}_\Omega).
\]
This and (\ref{CD:full}) imply that the spectral properties of $\LL$ 
on $\BB^{p,q}(T,V)$, $\KK$ on $\CC^{p,q}_\Omega$ and $\MM$ on  $\BBf_Z^{p,q}$
are (almost) identical. More precisely, the essential spectral radii and the eigenvalues 
of modulus larger than the essential spectral radius coincide, including multiplicity,
with an isometric bijection between the generalised eigenspaces. 

Recalling Subsection \ref{defs}, we decompose the operator $M=(\MM^m_{\vomega})_{\omega\omega'}$ into
\[
M_b=((\MM^m_{\vomega})_{\omega\omega'})_b
\quad\mbox{and}\quad 
M_c=((\MM^m_{\vomega})_{\omega\omega'})_c:\BBf^{\Theta_\omega,p,q}_\Gamma\to \hBBf^{\Theta_{\omega'},p,q}_\Gamma.
\]
From Proposition \ref{powerdecay}, the operator $((\MM^m_{\vomega})_{\omega\omega'})_c$ is compact. From Lemma \ref{lm:Mb}, it follows 
\[
\|((\MM^m_{\vomega})_{\omega\omega'})_b\|_{L(\BBf^{\Theta_\omega,p,q}_\Gamma,\hBBf^{\Theta_{\omega'} ,p,q}_\Gamma)}\le C \cdot 
\frac{\sup_{V_{\vomega}} |g^{(m)}|}
{\inf_{V_{\vomega}} \det(DT^m|_{E^u})}
\cdot \sup_{V_{\vomega}} \lambda^{(p,q,m)}\;.
\]
We decompose 
$\MM^m_{\vomega}:\BBf^{p,q}_Z\to \hBBf^{p,q}_Z$ into $
(\MM^m_{\vomega})_b, (\MM^m_{\vomega})_c:\BBf^{p,q}_Z\to \hBBf^{p,q}_Z$, by setting
\[
(\MM^m_{\vomega})_b((\mathbf{u}_\omega)_{\omega\in \Omega})
= \left(
\sum_{\omega\in \Omega} ((\MM^m_{\vomega})_{\omega\omega'})_b\mathbf{u}_\omega
\right)_{\omega'\in \Omega}\;
\]
and similarly for $(\MM^m_{\vomega})_c$. Finally we put
\[
(\MM^m)_b=\sum_{\vomega}(\MM^m_\vomega)_b \quad\mbox{and}\quad (\MM^m)_c=\sum_{\vomega}(\MM^m_\vomega)_c\; :\BBf^{p,q}_Z\to \hBBf^{p,q}_Z
\]
so that $\MM^m=(\MM^m)_b+(\MM^m)_c$. Then the operator $(\MM^m)_c$ is compact and
\begin{equation}\label{eqn:MMmb}
\|(\MM^m)_b\|_{L(\BBf^{p,q}_Z,\hBBf^{p,q}_Z)}
\le C \cdot \sum_{\vomega\in \Omega}
\frac{\sup_{V_{\vomega}} |g^{(m)}|}
{\inf_{V_{\vomega}} \det(DT^m|_{E^u})}
\cdot \sup_{V_{\vomega}} \lambda^{(p,q,m)}\;.
\end{equation}

\subsection{The end of the proof of Theorem~\ref{mainprop}}
\label{proofmainprop}
Since the spectral properties of $\LL$ 
on $\BB^{p,q}(T,V)$ and $\MM$ on  $\BBf_Z^{p,q}$
are (almost) identical as we noted, it is enough for the proof of Theorem~\ref{mainprop} to show that the essential spectral radius of $\MM$ on $\BBf_Z^{p,q}$ is bounded by $Q^{p,q}(T,g)=Q_{*}^{p,q}(T,g)$. Recall the positive-valued $C^{r-1}$ function $\tilde{g}$ taken just before the definition of the subsets $\Omega_m$. From standard argument in hyperbolic dynamical systems, there exists a constant $C(T,\tilde{g})>0$ such that
\begin{equation}\label{bdk}
\frac{\sup_{V_{\vomega}}\, \tilde{g}^{(m)}}{\inf_{V_{\vomega}}|\det(DT^m|_{E^u})|}
\le
C(T,\tilde{g})\cdot \inf_{V_{\vomega}} 
\left(\frac{\tilde{g}^{(m)}}{|\det(DT^m|_{E^u})|}
\right)
\end{equation}
for all $\vomega\in \Omega_m$ and $m\ge 1$. 
It follows
\[
\frac{\sup_{V_{\vomega}} g^{(m)}}
{\inf_{V_{\vomega}} \det(DT^m|_{E^u})}
\cdot \sup_{V_{\vomega}} \lambda^{(p,q,m)}
\le  C(T,\tilde{g})\cdot \sup_{V_{\vomega}} 
\left(\frac{\tilde{g}^{(m)}\cdot \lambda^{(p,q,m)}}{|\det(DT^m|_{E^u})|}
\right).
\]
Therefore we have, from (\ref{eqn:MMmb}), 
\begin{equation}\label{eqn:MbQ}
\|(\MM^m)_b\|_{L(\BBf^{p,q}_Z,\hBBf^{p,q}_Z)}
\le  C(T,\tilde{g}) \cdot 
Q^{p,q}_*(T,\tilde{g},\VV,m)\;.
\end{equation}
Since $(\MM^m)_c$ is compact, the essential spectral radius of $\MM:\BBf^{p,q}_Z\to \BBf^{p,q}_Z$ is bounded by $(C(T,\tilde{g}) \cdot 
Q^{p,q}_*(T,\tilde{g},\VV,m))^{1/m}$ and hence by $Q^{p,q}_*(T,\tilde{g})$, letting $m\to \infty$. 
This holds for any $C^{r-1}$ function $\tilde{g}$ such that $\tilde{g}(x)>|g(x)|$ on $X$. Therefore, by Lemma~\ref{lm:liminf},  the essential spectral radius of $\MM$ is bounded by  $Q^{p,q}_*(T,{g})$.

\begin{remark} 
We took  a positive-valued function $\tilde{g}$ (instead of $|g|$) so that (\ref{bdk}) holds. See Remark \ref{re:Ki12} also.  
\end{remark}

\section{The flat trace}\label{sec:trace}

In this section, we discuss about a flat trace for operators 
$\PP:\BBf^{p,q}_Z\to \hBBf^{p,q}_Z$ and give some related results. 
 
\subsection{Definition of the flat trace}
\label{flatrc}
Set $Z=\Omega \times \Gamma$. For $\zeta=(\omega,n,\sigma)\in Z$,   write
\[
n(\zeta)=n\, ,\quad \sigma(\zeta)=\sigma\, ,\quad \omega(\zeta)=\omega\quad\mbox{and}\quad \FF(\zeta)=\FF(\Theta_{\omega(\zeta)}), \quad \Theta(\zeta)=\Theta_{\omega(\zeta)}.
\] 
Then the Banach space $\BBf_Z^{p,q}$ introduced in the last section is written as 
\[
\BBf_Z^{p,q}:=\left\{ (u_{\zeta})_{\zeta\in Z}\;\left|\; 
u_{\zeta}\in \BBf^{\Theta(\zeta)}_{n(\zeta)+3}\;\mbox{ and }\; 
\lim_{n(\zeta)\to \infty} \left(2^{c(\sigma(\zeta))n(\zeta)} 
\|u_{\zeta}\|_{L^1(\FF(\zeta))} \right)=0 \right.
\right\}
\]
where $c(+)=p$ and $c(-)=q$. 
We will regard each element $\mathbf{u}$ of $\BBf_Z^{p,q}$ as a family $(u_{\zeta})_{\zeta\in Z}$ of functions with index set $Z$, and each $u_\zeta$ will be called the $\zeta$-component of $\mathbf{u}$.  
For $\zeta\in Z$, let 
$\BBf_\zeta$ (resp.$\hBBf_\zeta$) be the closed subspace of 
$\BBf^{p,q}_Z$ (resp.$\hBBf_Z$) that consists of elements $(u_{\zeta'})_{\zeta'\in Z}$ such that $u_{\zeta'}=0$ if $\zeta'\neq \zeta$, equipped with the restriction of the norm $\|\cdot \|_{\BBf^{p,q}_Z}$ (resp. $\|\cdot \|_{\hBBf^{p,q}_Z}$).

Consider a bounded operator $\PP:\BBf^{p,q}_{Z}\to \hBBf^{p,q}_{Z}$. For $\zeta,\zeta'\in Z$, let $
\PP_{\zeta\zeta'}: \BBf_{\zeta} \to \hBBf_{\zeta'}$ be the bounded operator that send $u\in \BBf_{\zeta}$ to the $\zeta'$-component of $\PP(u)$. Observe that the restriction of $\PP_{\zeta\zeta'}$ to $\hBBf_{\zeta}$ is written as an integral operator with kernel
\[
K_{\zeta\zeta'}(x,y)=\PP_{\zeta\zeta'} (\widehat \chi_{n(\zeta)+2} (y-\cdot)\bigr )
(x)\,.
\]
Indeed, for $u\in \hBBf_\zeta$, 
we have
\[
\PP_{\zeta\zeta'} u(x)=\PP_{\zeta\zeta'} (\widehat{\chi}_{n(\zeta)+2}*u)(x)=\int \PP_{\zeta\zeta'} (\widehat \chi_{n(\zeta)+2} (\cdot-y))(x) \cdot  u(y) dy. 
\]
Since $\widehat \chi_{n(\zeta)+2} (\cdot-y)$ belongs to $\BBf_\zeta$ and depends on $y\in \real^d$ continuously, the kernel $K_{\zeta\zeta'}(x,y)$  
is continous with respect to $x$ and $y$. If 
$K_{\zeta\zeta}(x,x)$ is integrable with respect to $x$,
we say that $\PP_{\zeta\zeta}$ admits a flat trace and put
\[
\tr^\flat(\PP_{\zeta\zeta})=\int_{\real^d} K_{\zeta\zeta}(x,x) dx.
\] 
\begin{remark}
The operator $\PP_{\zeta\zeta}$ may be expressed as  integral operators with different kernels. 
And the different choice of kernels may give different traces for $\PP_{\zeta\zeta}$. 
\end{remark}

\begin{definition}
We say that a bounded operator $\PP:\BBf^{p,q}_Z\to \hBBf^{p,q}_Z$ admits a flat trace if 
$\PP_{\zeta\zeta}$ for each $\zeta\in Z$ admits a flat trace and if the following limit exists:
\[
\tr^\flat \PP:=\lim_{n\to \infty} \sum_{\zeta: n(\zeta)\le n} \tr^\flat (\PP_{\zeta\zeta}).
\]
If $\PP^m$ admits a flat trace for all $m\ge 1$,
we define the  flat
determinant of $\PP$ to be the formal power series
\begin{equation}
{\det}^\flat(\id-z\PP)
=\exp-\sum_{m\ge 1}\frac{z^m}{m}
{\tr}^\flat(\PP^m)\, .
\end{equation}
\end{definition}
Clearly, if  $\tr^\flat(\PP_1)$ and
$\tr^\flat(\PP_2)$ are well-defined, then so is $\tr^\flat(\PP_1+\PP_2)$, and
$\tr^\flat (\PP_1)+\tr^\flat (\PP_2)=\tr^\flat(\PP_1+\PP_2)$.
\begin{proposition}\label{nucltrace}
Suppose that
$\PP:\BBf^{p,q}_{Z}\to \hBBf^{p,q}_{Z}$ is a bounded operator and 
has a nuclear representation 
$\PP=\sum_{i=1}^{\infty} v_i\otimes v_i^*$ where $v_i\in \hBBf^{p,q}_{Z}$ and 
$v_i^*\in (\BBf^{p,q}_{Z})^*$ 
satisfy $\sum_{i=1}^{\infty} \|v_i\|_{\hBBf^{p,q}_{Z}}\cdot \|v_i^*\|_{\BBf^{p,q}_{Z}}<\infty$. 
Then $\PP$ admits a flat trace. Further it holds
\[
\tr^\flat \PP= \sum_{i=1}^{\infty} v_i^*(v_i)\, \quad\mbox{and}\quad
|\tr^\flat \PP|\le  \sum_{i=1}^{\infty} \|v_i\|_{\hBBf^{p,q}_{Z}}\cdot \|v_i^*\|_{\BBf^{p,q}_{Z}}<\infty\,.
\]
\end{proposition}
\begin{proof} Put $\PP^{(i)}=v_i\otimes v_i^*$. 
Let $v_{i,\zeta}$ for $\zeta\in Z$ be the $\zeta$-component of $v_i$. Also let $v^*_{i,\zeta}$ be the functional 
on $\BBf_\zeta$ that $v_i^*$ induces. Then we have $\PP^{(i)}_{\zeta\zeta}=v_{i,\zeta}\otimes v^*_{i,\zeta}$ and also
\footnote{For the second equality, recall Remark \ref{rem:B}. This equality would not hold if we used the Banach space ${\tilde{\BBf}}^{\Theta,p,q}_\Gamma$ in the place of ${\BBf}^{\Theta,p,q}_\Gamma$ } 
\begin{equation}\label{eqn:supv}
\sup_{\zeta\in Z} \|v_{i,\zeta}\|_{\BBf^{p,q}_{Z}}=\|v_i\|_{\BBf^{p,q}_{Z}}\quad\mbox{and}\quad
\sum_{\zeta\in Z} \|v_{i,\zeta}^*\|_{\BBf^{p,q}_{Z}}= \|v^*_i\|_{\BBf^{p,q}_{Z}}. 
\end{equation}
By definition we have
\[
\tr^\flat \PP^{(i)}_{\zeta\zeta}=\int K^{(i)}_{\zeta\zeta}
(x,x) \, dy\quad\mbox{where}\quad K^{(i)}_{\zeta\zeta}(x,y)=v_{i,\zeta}(x) v^*_{i,\zeta}(\widehat \chi_{n(\zeta)+2} (\cdot-y)).
\]
Since $\widehat \chi_{n(\zeta)+2} (\cdot-y)$ for $y\in \real^d$ is uniformly bounded in $\BBf_{\zeta}$, we have, by (\ref{eqn:cr2}), that  
\[
\int |K^{(i)}_{\zeta\zeta}(x,x)| dx \le C(\zeta)\cdot  \|v_i\|_{\hBBf^{p,q}_{Z}}\cdot \|v_i^*\|_{\BBf^{p,q}_{Z}}\;\quad \mbox{for all $1\le i<\infty$.}
\]
This implies that $\sum_i K^{(i)}_{\zeta\zeta}(x,x)$ is integrable with respect to $x$, that is,  $\PP_{\zeta\zeta}$ admits a flat trace. 
Since  $\widehat{\chi}_{n(\zeta)+2}* v_{i,\zeta}=v_{i,\zeta}$ for $v_{i,\zeta}\in \BBf^{\Theta(\zeta)}_{n(\zeta)+2}$,  it holds
\[
\tr^\flat \PP^{(i)}_{\zeta\zeta}=v^*_{i,\zeta}\left(\int v_{i,\zeta}(x)\widehat \chi_{n(\zeta)+2} (\cdot-x) \, dx
\right)=v_{i,\zeta}^*(v_{i,\zeta})
\]
for each $\zeta\in Z$ and $i\ge 1$.
It follows from  (\ref{eqn:supv}) that 
\[
\sum_i\sum_{\zeta\in Z}| \tr^\flat \PP^{(i)}_{\zeta\zeta}| \le
\sum_i\sum_{\zeta\in Z} \|v_{i,\zeta}\|_{\hBBf^{p,q}_{Z}}
\|v_{i,\zeta}^*\|_{\BBf^{p,q}_{Z}}
\le \sum_i\|v_i^*\|_{\hBBf^{p,q}_{Z}} \|v_i\|_{\BBf^{p,q}_{Z}}<\infty\,.
\]
Therefore we conclude that $\tr^\flat \PP$ exists and 
\[
\tr^\flat \PP=\lim_{n\to \infty} \sum_{n(\zeta)\le n}\sum_i\tr^\flat \PP^{(i)}_{\zeta\zeta} =\sum_i \sum_{\zeta\in Z}v_{i,\zeta}^*(v_{i,\zeta})
=\sum_i v_{i}^*(v_{i})\,.
\]
The inequality for $|\tr^\flat \PP|$ is then obvious. 
\end{proof}

\subsection{The flat trace of the operators $\MM^m$}
We next consider the flat traces of the operators $\MM^m$, $(\MM^m)_b$ and $(\MM^m)_c$ introduced in the last section. 
The flat trace of $\MM^m$ coincides with the dynamical trace:
\begin{proposition}\label{dyntrace}
The operator $\MM^m:\BBf^{p,q}_Z\to \hBBf^{p,q}_Z$ 
for $m\ge 1$ admits a flat trace and holds 
\[
\tr^\flat(\MM^m)
=\sum_{T^m (x)=x} \frac{g^{(m)}(x)}
{|\det(\id -DT^m(x))|}\, .
\]
\end{proposition}

\begin{proof}
Consider $\MM^m$ for $m\ge 1$. Take $\omega\in \Omega$ and $\vomega\in \Omega_m$ and  recall the definition of the operator  $(\MM^m_{\vomega})_{\omega\omega}$. 
Then we see that, for each $\zeta\in Z$ with $\omega(\zeta)=\omega$, the flat trace $\tr^\flat(\MM^m_{\zeta\zeta})$ is defined as the integral
\[
\int \widehat{\psi}_{\Theta(\zeta),n(\zeta),\sigma(\zeta)}(x-y) \cdot G(y) \cdot \widehat{\tilde{\psi}}_{\Theta(\zeta),n(\zeta),\sigma(\zeta)}(\TT(y)-z)\cdot \widehat{\chi}_{n(\zeta)+2}(z-x) dx dy dz
\]
where $\TT$ and $G$ are those in the setting (\ref{setting}) with $\omega'=\omega$.
Since  
\[
\widehat{\tilde{\psi}}_{\Theta(\zeta),n(\zeta),\sigma(\zeta)}*\widehat{\chi}_{n(\zeta)+2}*\widehat{\psi}_{\Theta(\zeta),n(\zeta),\sigma(\zeta)}=\widehat{\psi}_{\Theta(\zeta),n(\zeta),\sigma(\zeta)}
\]
by (\ref{eqn:psed}), 
we see that $\MM^m_{\zeta\zeta}$ admits a flat trace and 
\[
\tr^\flat(\MM^m_{\zeta\zeta})=\int \widehat{\psi}_{\Theta(\zeta),n(\zeta),\sigma(\zeta)}(\TT(x)-x) \cdot G(x)  \, dx\,.
\]
Thus, for each integer $n_0$, we have 
\begin{equation}\label{eqn:trlocal}
\sum_{\zeta: n(\zeta)\le n_0; \omega(\zeta)=\omega}
\tr^\flat(\MM^m_{\zeta\zeta})
=\int \widehat{\chi}_{n_0}(\TT(x)-x) \cdot G(x)  \, dx\, .
\end{equation} 
The function $\widehat{\chi}_{n_0}$,  regarded as a distribution, converges to the Dirac measure at $0$ as $n_0\to \infty$. 
Note that there is at most one 
fixed point of $\TT$ in $\supp(G)$ because the covering $\VV$ is assumed to be generating. If there is no fixed point in $\supp(G)$, the sum (\ref{eqn:trlocal}) converges to zero as $n_0\to \infty$. If there is one fixed point $x_0$ in $\supp(G)$, that fixed point should be hyperbolic by hyperbolicity of $T$ and hence we may perform a local change of variable $z=\TT(x)-x$ in its small neighborhood to obtain 
\[
\lim_{n_0\to \infty} \int \widehat{\chi}_{n_0}(\TT(x)-x) \cdot G(x)  dx =\frac{G(x_0)}
{|\det(\id -D\TT(x_0))|}\,.
\]
Recalling the definition of $\TT$, $G$ and $h_\omega$, we see that the operator $\MM^m_{\vomega}$ admits a flat trace and 
\[
\tr^\flat(\MM^m_{\vomega})=\sum_{\omega\in \Omega} \sum_{T^m(x)=x} \frac{\phi_{\omega}(x)\cdot \phi_{\vomega}(x)\cdot g^{(m)}(x)}
{|\det(\id -DT^m(x))|}=
\sum_{T^m(x)=x} \frac{\phi_{\vomega}(x)\cdot g^{(m)}(x)}
{|\det(\id -DT^m(x))|}.
\]
Taking the sum with respect to $\vomega\in \Omega_m$, we obtain the proposition. \end{proof}

The following property of $(\MM^m)_b$ is important in the proof of Theorem~\ref{main}.
\begin{proposition}\label{zerotrace}
There is 
 $L=L(T,g)\ge 1$  so that,
if $m_j\ge L$ for $1\le j \le J$, then
\[
\biggl(\prod_{j=1}^J(\MM^{m_j})_b \biggr)_{\zeta\zeta}=0 
\quad \mbox{ for all $\zeta\in Z$,}
\]
and in particular 
\[
\tr^\flat \biggl ( \prod_{j=1}^J(\MM^{m_j})_b \biggr ) =0\, 
\,.
\]
\end{proposition}
\begin{remark}\label{rem:prod}
We read the expression $\Pi_{i=1}^{k-1} \PP_i$ as the product $\PP_{k-1}\PP_{k-2}\cdots \PP_1$, not as $\PP_{1}\PP_{2}\cdots \PP_{k-1}$.
\end{remark}
\begin{proof}
By hyperbolicity of $T$, there exists $L\ge 1$, such that, for all $\omega, \omega'\in \Omega$ and $\vomega\in \Omega_m$ with $m\ge L$, we have 
$h_{\max}^+(\TT,G)<0$ and $h_{\min}^{-}(\TT,G) >0$ in the setting (\ref{setting}). 
This and the definition of the relation $\hookrightarrow$ imply that $(\MM^m)_b$ for $m\ge L$ is "strictly lower triangular" as a matrix of operators in the sense that $((\MM^{m})_b)_{\zeta\zeta'}\neq 0$ only if $\sigma(\zeta')n(\zeta')<\sigma(\zeta)n(\zeta)$. 
Clearly this gives the claim of the lemma. 
  \end{proof}

For the operator $(\MM^m)_c$, we have the following: Recall Subsection \ref{appr} and put
\begin{equation}\label{K}
k_*=k_*(d,r,p,q):=\left  [1+\frac{d}{r-p+q-1} \right ] \, .
\end{equation} 
\begin{proposition}\label{Knucl}
The operator $(\MM^m)_c$ belongs to the class $\LL^{(a)}_t(\BBf^{p,q}_{Z}, \hBBf^{p,q}_{Z})$ for any $t>d/(r-p+q-1)$. Further, 
given any bounded operators $\PP_j:\BBf^{p,q}_Z\to \hBBf^{p,q}_Z$ for $0\le j\le k_*$
and any integers $m_j\ge 1$ for $1\le j\le k_*$, we have  
\[
\biggl\|
\PP_0 \prod_{j=1}^{k_*}((\MM^{m_j})_c\PP_j) \biggr\|_{1}^{(a)}
\le C  \cdot 
\|\PP_0\|_{L(\BBf^{p,q}_{Z}, \hBBf^{p,q}_{Z})}  
\prod_{j=1}^{k_*} \left ( \|(\MM^{m_j})_c\|_{k_*}^{(a)} \cdot \|\PP_j\|_{L(\BBf^{p,q}_{Z}, \hBBf^{p,q}_{Z})}\right )
\]
\end{proposition}
\begin{proof}
The first claim is a consequence of Proposition \ref{powerdecay} and the definition of the operator $(\MM^{m})_c$. The second then follows from (\ref{prod2}) and (\ref{mult}). 
\end{proof}

\section{Dynamical determinants: Proof of Theorem~\ref{main}}
\label{S5}
In this section, we prove Theorem~\ref{main}. 
Let $q<0<p$ be so
that $p-q<r-1$.
As we noted in Subsection~\ref{sec:MM}, the operators $\MM$ on $\BBf^{p,q}_Z$ and
$\LL$ on $\BB^{p,q}(T,V)$ share almost same spectral data.  
And we have proved in Subsection~\ref{proofmainprop} that the essential spectral
radius of $\MM$ is not larger than
$Q_*^{p,q}(T,g)=Q^{p,q}(T,g)$. Fix $\epsilon>0$ arbitrarily and set $\rho:=Q^{p,q}(T,g)+\epsilon$. 
By Proposition~\ref{dyntrace} we have $d_\LL(z)={\det}^\flat (\id-z\MM)$  
as formal power series.
Therefore,
in order to prove Theorem~\ref{main}, it suffices to show that
$d_\LL(z)={\det}^\flat (\id-z\MM)$ extends holomorphically to the
disc of radius $\rho^{-1}$,  and that $d_\LL(z)$ vanishes at
order $n_z$ in this disc if and only if $1/z$ is an eigenvalue of
algebraic multiplicity $n_z$ for
$\MM$ on $\BBf^{p,q}_Z$.
This is the content of the present section.

Consider the spectral projector $\PP_0$ for $\MM:\BBf^{p,q}_Z\to \BBf^{p,q}_Z$
associated to eigenvalues of modulus larger than or equal to
$\rho$.
We have\footnote{Since $(\MM_0)^m=\MM^m\PP_0=(\MM^m)_0$, we
write $\MM^m_0$ without risk of confusion, similarly for $\MM_1$.}
$\MM^m=\MM_0^m+\MM_1^m$, with $\MM_0=\MM\PP_0$
and $\MM_1=\MM(\id-\PP_0)$. 
For each $m\ge 1$, the operator $\MM_0^m$ is of finite rank and its image is contained
 in $\hBBf^{p,q}_Z$. By Proposition~\ref{nucltrace},  $\MM_0^m$ admits a flat trace and its flat trace coincides with the usual trace defined for finite rank operators. 
By Proposition~\ref{dyntrace}, we also see that $\MM_1^m=\MM^m-\MM_0^m$ admits a flat trace.  Hence we 
may decompose
\[
{\det}^\flat(\id-z\MM)={\det}^\flat(\id-z\MM_0){\det}^\flat(\id-z\MM_1)\, 
\]
and the factor 
$
{\det}^\flat(\id-z\MM_0) 
$
is  a polynomial which vanishes exactly
at the inverse eigenvalues of $\MM\PP_0$, with order equal to the multiplicity
of the eigenvalue. 
To
prove Theorem~\ref{main}, it thus
suffices to show that $\det^\flat(\id-z\MM_1)$ is
holomorphic and nowhere zero in the disc of radius $\rho^{-1}$, i.e., 
for any $\epsilon'>0$, there exists $C>0$ such that 
\begin{equation}\label{trbd}
|\tr^{\flat}(\MM_1^m)|<C (\rho+\epsilon')^{m}\quad \mbox{for all $m
\ge 1$.}
\end{equation}

Since the proof is much simpler in the case $r>d+p-q+1$, we will discuss about such case first in 
Subsections~\ref{largerr} and consider the other case later in Subsection~\ref{smallr}.  
From (\ref{eqn:MbQ}), we may take an integer $m_0\ge 1$ so that 
\begin{equation}\label{mzero}
\|(\MM^m)_b\|_{L(\BBf^{p,q}_Z,\hBBf^{p,q}_Z)}\le \rho^m\quad\mbox{ and }\quad
\|\MM_1^m\|_{L(\BBf^{p,q}_Z,\hBBf^{p,q}_Z)}\le \rho^m\quad \mbox{ for $m\ge m_0$. }
\end{equation}
For $m\ge 1$, we put $
(\MM^m)_a:=\MM_1^m-(\MM^m)_b=(\MM^m)_c-\MM_0^m$
 so that 
\begin{equation}\label{eqn:abdec}
\MM_1^m=(\MM^m)_a+ (\MM^m)_b.
\end{equation} 

\subsection{ The case $r>p-q+d+1$}
\label{largerr}
By Proposition \ref{Knucl}, the operators $(\MM^m)_c$ and $(\MM^m)_a$ both belong to
$\LL^{(a)}_1(\BBf^{p,q}_Z,\hBBf^{p,q}_Z)$ in this case. Write 
a large  integer $n$ as a sum $n=m(1)+m(2)+\dots+m(k)$ with $m_0\le m(i)\le 2m_0$. Then, using (\ref{eqn:abdec}), 
we may write
 $\MM_1^n=\prod_{i=1}^{k}\MM_1^{m(i)}$ as
\[
\MM_1^n=\prod_{i=1}^{k}(\MM^{m(i)})_b+\sum_{j=1}^{k} 
\left(\prod_{i=j+1}^{k}\MM_1^{m(i)} \cdot (\MM^{m(j)})_a \cdot \prod_{i=1}^{j-1}(\MM^{m(i)})_b
\right)\, .
\]
By Proposition~\ref{zerotrace}, we have 
$\tr^{\flat}\left( \prod_{i=1}^{k}(\MM^{m(i)})_b \right)=0$ for the first term.
For the other terms, we have, by Proposition~\ref{Knucl} with $k_*=1$ and (\ref{mzero}), that 
\begin{align*}
&\left|\tr^{\flat} \left( \prod_{i=j+1}^{k}\MM_1^{m(i)}\cdot (\MM^{m(j)})_a 
\cdot \prod_{i=1}^{j-1}(\MM^{m(i)})_b
\right)\right|\\
&\qquad\le \prod_{i=j+1}^{k}\|\MM_1^{m(i)}\|_{\BBf^{p,q}_Z}\cdot 
\|(\MM^{m(j)})_a\|^{(a)}_{1} 
\cdot \prod_{i=1}^{j-1}\|(\MM^{m(i)})_b\|_{\BBf^{p,q}_Z} 
\le C(m_0) \rho^{n} \, .
\end{align*}
Therefore we obtain the claim (\ref{trbd}).


\subsection{The case $r\le d+p-q+1$}
\label{smallr}
In this case, the operators $(\MM^m)_c$ and $(\MM^m)_a$ may not belong to
$\LL^{(a)}_1(\BBf^{p,q}_Z,\hBBf^{p,q}_Z)$, so we have to modify the simple proof in the 
last subsection.  Recall the integer $k_*=k_*(d,r,p,q)\ge 2$ from (\ref{K}).

Consider a large integer $n$ and write it as a sum $n=m(1)+m(2)+\dots+m(k)$ with $m_0\le m(i)\le 2m_0$. Using (\ref{eqn:abdec}), 
we write the product $\MM_1^n=\prod_{i=1}^{k}\MM_1^{m(i)}$ as
\begin{align}\label{MMn}
\MM_1^{n}=\prod_{i=1}^{k}(\MM^{m(i)})_b\,\,
&+\sum_{\nu< k_*}\;\;\sum_{1\le j(1)<j(2)<\dots<j(\nu)\le k} M(\{j(\ell)\}_{\ell=1}^{\nu})\\
&+\;\;
\sum_{1\le j(1)<j(2)<\dots<j(k_*)\le k} M'(\{j(\ell)\}_{\ell=1}^{k_*})\, ,\notag
\end{align}
where, setting $j(0)=0$,  
\[
M(\{j(\ell)\}_{\ell=1}^{\nu})=\prod_{i=j(\nu)+1}^{k} (\MM^{m(i)})_b\cdot \prod_{\ell=1}^{\nu}\biggl(
(\MM^{m(j(\ell))})_a \cdot
\biggl(
\prod_{i=j(\ell-1)+1}^{j(\ell)-1}(\MM^{m(i)})_b 
\biggr)
\biggr)\, ,
\]
and
\[
M'(\{j(\ell)\}_{\ell=1}^{k_*})=\prod_{i=j(k_*)+1}^{k} \MM^{m(i)}_1\cdot \prod_{\ell=1}^{k_*}
\biggl((\MM^{m(j(\ell))})_a\cdot 
\biggl(
\prod_{i=j(\ell-1)+1}^{j(\ell)-1}(\MM^{m(i)})_b 
\biggr)
\biggr) \, .
\]
\begin{remark}
The decomposition above is obtained as follows. Consider the process to expand $
\MM_1^n=\prod_{i=1}^{k}\MM_1^{m(i)}$ 
 using  (\ref{eqn:abdec}) for $m=m(i)$  in the turn $i=1,2,, \dots, k$. For instance, we have 
\[
\MM_1^n= \left(\prod_{i=2}^{k}\MM_1^{m(i)}\right)\cdot (\MM^{m(1)})_a
+\left(\prod_{i=2}^{k}\MM_1^{m(i)}\right)\cdot (\MM^{m(1)})_b
\]
for the first step. When we find a term that contains $(\MM_1^{m(i)})_a$ for $k_*$ times, proceeding in this way, we stop expanding that term, obtaining the terms in the second sum in (\ref{MMn}). The other resulting terms are collected in the first sum.
\end{remark}

From Proposition \ref{zerotrace}, the flat trace of the first term on the right hand side of  (\ref{MMn}) is zero. Therefore,  to prove (\ref{trbd}), it suffices to show the following estimates for the other terms:
\begin{equation}\label{eqn:Mj1}
|\tr^{\flat} M'(\{j(\ell)\}_{\ell=1}^{k_*})|
\le C \rho^{n}
\end{equation}
and
\begin{equation}\label{eqn:Mj2}
|\tr^{\flat} M(\{j(\ell)\}_{\ell=1}^{\nu})|
\le C \rho^{n}\,.
\end{equation} 
By Proposition~\ref{Knucl} and (\ref{mzero}), we can see that  $M'(\{j(\ell)\}_{\ell=1}^{k_*})\in \LL^{(a)}_1(\BBf^{p,q}_Z,\hBBf^{p,q}_Z)$ and that the estimate (\ref{eqn:Mj1}) holds.  
Since $(\MM^{m(i)})_c-(\MM^{m(i)})_a=\MM_0^{m(i)}$ is of finite rank, the estimate (\ref{eqn:Mj2}) follows if we show
\begin{equation}\label{keym}
\biggl |\tr^{\flat}\biggl(\prod_{i=j(\nu)+1}^{k} (\MM^{m(i)})_b\cdot \prod_{\ell=1}^{\nu}
\biggl( (\MM^{m(j(\ell))})_c  \cdot 
\prod_{i=j(\ell-1)+1}^{j(\ell)-1}(\MM^{m(i)})_b 
\biggr)
\biggr)\biggr |\le 
C \rho^{n}\,.
\end{equation}

In the following, we will work directly with kernels of  operators to prove (\ref{keym}). 
Although the notation become a little complex, the argument is straightforward. 
Let $\mathcal{Y}$ be the set of sequences   
$\{\omega(i), \vomega(i)\}_{i=1}^{k}$ with $\omega(i)\in \Omega$ and 
$\vomega(i) \in \Omega_{m(i)}$. 
For $Y=\{\omega(i), \vomega(i)\}_{i=1}^{k}\in \mathcal{Y}$ and $1\le i\le k$, we define the relation $\hookrightarrow_{Y,i}$ on $\Gamma$ as the relation $\hookrightarrow_{\TT,G}$ defined in the setting (\ref{setting}) with $\omega=\omega(i)$, $\omega'=\omega(i+1)$ and $\vomega=\vomega(i)$. (Put $\omega(k+1)=\omega(1)$ in case $i=k$.)  
For $Y=\{\omega(i), \vomega(i)\}_{i=1}^{k}\in \mathcal{Y}$,
let $\mathcal{Z}(Y)$ be the set of sequences $\{\zeta(i)\}_{i=1}^{k+1}$ in $Z^{k+1}$ such that $\zeta(k+1)=\zeta(1)$, that $\omega(\zeta(i))=\omega(i)$ for all $1\le i\le k$ and that 
\[
(n(\zeta(i)),\sigma(\zeta(i))) \hookrightarrow_{Y,i} (n(\zeta(i+1)),\sigma(\zeta(i+1)))\quad \mbox{iff}\quad  i\notin J:=\{j(\ell)\}_{\ell=1}^{\nu}\,.
\]  
Then, to prove (\ref{keym}), it is enough to show 
\begin{equation}\label{eqn:trsum}
\sum_{Y=\{\omega(i),\vomega(i)\}_{i=1}^{k}\in \mathcal{Y}}\;\;
\sum_{\{\zeta(i)\}_{i=1}^{k+1}\in \mathcal{Z}(Y)}\;\;
\left|\tr^{\flat}
\left(
\prod_{i=1}^{k}
(\MM^{m(i)}_{\vomega(i)})_{\zeta(i)\zeta(i+1)}
\right)\right|\le C\rho^n\,. 
\end{equation}

Let  $Y=\{\omega(i), \vomega(i)\}_{i=1}^{k}\in \mathcal{Y}$ and $\{\zeta(i)=(\omega(i),n(i),\sigma(i))\}_{i=1}^{k+1}\in \mathcal{Z}(Y)$. 
By the definition, we have
\[
\tr^{\flat} \left( \prod_{i=1}^{k} 
(\MM_{\vomega(i)}^{m(i)})_{\zeta(i)\zeta(i+1)}\right)= \int \left(\prod_{i=1}^{k} K_{i}(x_i,x_{i+1})\right) dx_1 dx_2 \dots dx_{k}
\]
where we read $x_{k+1}=x_1$ and put
\[
K_{i}(x,y)=\int \widehat{\psi}_{\Theta_{\omega(i+1)}, n(i+1), \sigma(i+1)}(y-w)\cdot  G_i(w)\cdot  \widehat{\tilde{\psi}}_{\Theta_{\omega(i)}, n(i), \sigma(i)}(\TT_i(w)-x) dw 
\]
with  $\TT_i:=\TT^{m(i)}_{\vomega(i), \omega(i+1)\omega(i)}$ and $G_i:=G_{\vomega(i), \omega(i+1)\omega(i)}$. (Recall Subsection \ref{system}.)
Here we canceled  the term $\widehat{\chi}_{n(1)+2}$ by using $\widehat{\tilde{\psi}}_{\Theta_{\omega(1)}, n(1), \sigma(1)}*\widehat{\chi}_{n(1)+2}=\widehat{\tilde{\psi}}_{\Theta_{\omega(1)}, n(1), \sigma(1)}$.

For all $1\le i\le k$, we have, from (\ref{eqn:kerpsi}), 
\[
|K_{i}(x,y)|\le C\cdot \|G_i\|_{L^\infty}\cdot \int b_{n(i+1)}(y-w) \cdot b_{n(i)}(\TT_i(w)-x) \, dw\,.
\] 
If $i\in J$, we have, from (\ref{inte}), 
\[
|K_{i}(x,y)|\le 
C(\TT_i,G_i) \cdot 2^{-(r-1)\max\{n(\zeta(i)), n(\zeta(i+1))\}} b_{\min\{n(\zeta(i)), n(\zeta(i+1))\}}(y-x).
\]
Therefore, using (\ref{eqn:bconv}), we see that $|\tr^{\flat} ( \prod_{i=1}^{k} 
(\MM_{\vomega(i)}^{m(i)})_{\zeta(i)\zeta(i+1)})|$ is bounded by 
\begin{align*}
&C^k\cdot C(m_0)^{k_*}\cdot   \prod_{i\in J}2^{-(r-1)\max\{n(i+1),n(i)\}} \prod_{i\notin J}
\sup_{V_{\vomega(i)}} |g^{(m(i))}|\\
&
\cdot \int {b}_{\nu(k)}(\TT_{k}(x_{0})-x_k) 
\cdots 
b_{\nu(2)}(\TT_2(x_3)-x_2)
\cdot 
b_{\nu(1)}(\TT_1(x_2)-x_1)
\, dx_1 \cdots dx_{k}\, ,
\end{align*}
where the constant $C$ does not depend on the choice of $m_0$ while $C(m_0)$ may, and  
\[
\nu(i)=\begin{cases}\min\{n(i),n(i+1)\},&\mbox{ if $i\in J$ };\\
n(i), &\mbox{ if $i\notin J$ }.
\end{cases}
\]  
From hyperbolicity of $T$ and the assumption that the covering $\VV$ is generating, we may choose the extensions of the diffeomorphisms $\TT_i$ so that the mapping 
\[
S:(\real^d)^k\to (\real^d)^k, \quad (x_i)_{i=1}^{k}\mapsto (\TT_1(x_2)-x_{1},\TT_2(x_3)-x_{2},\dots,\TT_{k}(x_{0})-x_{k})\, ,
\]
is a diffeomorphism and satisfy 
\[
\inf_{\real^d} |\det DS|\ge C^{-k}\cdot  \prod_{i=1}^{k}\left(\sup_{V_{\vomega(i)}}  |\det DT^{m(i)}|_{E^u}|\right)
\]
for a constant $C>0$ that does not depend on the choice of $m_0$. 
Hence we get the following estimate for each term in (\ref{eqn:trsum}):
\begin{align}\label{eqn:eachtr}
\biggl | \tr^{\flat}&\left( \prod_{i=1}^{k} (\MM_{\vomega(i)}^{m(i)})_{\zeta(i)\zeta(i+1)}\right)
\biggr |\\
&\le 
C^k\cdot C(m_0)^{k_*}\prod_{i\in J}2^{-(r-1)\max\{n(i+1),n(i)\}} 
\prod_{i\notin J} \frac
{\sup_{V_{\vomega(i)}}|g^{(m(i))}|}
{\inf_{V_{\vomega(i)}}|\det DT^{m(i)}|_{E^u}|}\,. \notag
\end{align}

Putting $c(-)=q$ and $c(+)=p$, we have,
for $i\notin J$,
\[
c(\sigma(i+1))n(i+1)-c(\sigma(i))n(i)\le 
\min\{p \cdot h_{\max}^{+}(\TT_i,G_i), q \cdot h_{\min}^{-}(\TT_i,G_i)\}\,.
\]
We may assume that the right hand side is negative, taking larger $m_0$ if necessary.
Since $\zeta(k+1)=\zeta(1)$, we have 
\begin{align}\label{eqn:cn2}
-&\sum_{i\in J}(r-1)\max\{n(i+1),n(i)\}\\
&\le -\left(\sum_{i\in J}(c(\sigma(i+1))n(i+1)-c(\sigma(i))n(i))\right)
+(r-p+q-1) \max_{1\le i\le k}n(i)
\notag
\\
&\le  \left(\sum_{i\notin J}(c(\sigma(i+1))n(i+1)-c(\sigma(i))n(i))\right)
+(r-p+q-1) \max_{1\le i\le k}n(i).
\notag
\end{align}
Therefore we obtain 
\begin{align*}
&
\biggl |\sum_{\{\zeta(i)\}_{i=1}^{k+1}\in \mathcal{Z}(Y)} 
\tr^{\flat}\left( \prod_{i=1}^{k} (\MM^{m(i)}_{\vomega(i)})_{\zeta(i)\zeta(i+1)}\right)\biggr | \\
&\le 
C^k \cdot C(m_0)^{k_*}
\cdot \prod_{i\notin J}
\biggl(
  \frac{\sup_{V_{\vomega(i)}}|g^{(m(i))}|}
  {\inf_{V_{\vomega(i)}}|\det DT^{m(i)}|_{E^u}|}\cdot 
2^{\min\{ p h_{\max}^{+}(\TT_i,G_i), q h_{\min}^{-}(\TT_i,G_i)\}}
\biggr)\\
&\le 
C^k \cdot C(m_0)^{k_*}
\cdot \prod_{i\notin J}
\biggl(
\frac{\sup_{V_{\vomega(i)}}|g^{(m(i))}|}
{\inf_{V_{\vomega(i)}}|\det DT^{m(i)}|_{E^u}|}\cdot 
\sup_{V_{\vomega(i)}} \lambda^{(p,q,m)}
\biggr).
\end{align*}
Recalling the function $\tilde{g}$ and the definition of $\Omega_m$ in Subsection \ref{localcharts} and using (\ref{bdk}), we conclude that (\ref{eqn:trsum}) is bounded by 
\[
C(\tilde{g})^k \cdot C(m_0)^{k_*} \cdot \prod_{i\notin J}(Q_*^{p,q}(T,\tilde{g},\mathcal{W},m(i)))\, .
\]
By Lemma \ref{lm:tp}, we may take $\tilde{g}$ so that $Q_*^{p,q}(T,\tilde{g})<\rho$. 
Since the constant $C(\tilde{g})$ above does not depend on the choice of $m_0$, we obtain (\ref{keym}) by taking large $m_0$.  We finished the proof of Theorem \ref{main}.

\appendix

\section{Eigenvalues and eigenvectors for different Banach spaces}\label{apd:eigenvalue}

In Theorem ~\ref{mainprop}, we may choose a variety of $p$ and $q$. Besides, as we will see in 
the proof, the space $\BB^{p,q}(T,V)$ depends on many objects, 
such as the system of local charts. Moreover, in \cite{BT} and \cite{GoL}, other 
Banach spaces of distribution were introduced, for which the analogue of 
Theorem~ \ref{mainprop} holds with  different bounds on the essential 
spectral radius. 
So one may ask to what extent the eigenvalues of the Ruelle transfer operator 
on different Banach spaces coincide. Theorem \ref{main} gives one answer to this question 
because the dynamical Fredholm determinant does not depend on the choice of Banach spaces. 
The following simple abstract lemma, which may not be new, gives a more direct answer.

\begin{lemma}\label{changespace}
Let $B$ be a separated topological linear space and let 
$(B_1,\|\cdot\|_1)$ and $(B_2,\|\cdot\|_2)$ be Banach spaces that are continuously embedded in $B$. Suppose further that
there is a subspace $B_0\subset B_1\cap B_2$ that is dense both in the Banach spaces 
$(B_1,\|\cdot \|_1)$ and $(B_2, \|\cdot \|_2)$. 
Let $L:B\to B$ be a continuous linear map, which preserves the subspaces $B_0$, $B_1$ and $B_2$. 
Suppose that the restriction of $L$ to $B_1$ and $B_2$ are bounded operators whose essential 
spectral radii are both strictly smaller than some number $\rho>0$.  Then the eigenvalues of 
$L:B_1\to B_1$ and $L:B_2\to B_2$ in $\{z\in \complex\mid |z|>\rho\}$ coincide. Further the 
corresponding generalized eigenspaces coincide and are contained in the intersection $B_1\cap B_2$. 
\end{lemma}

\begin{proof}
First, we show that the essential spectral radius $r_{ess}(\LL)$ of an operator 
$\LL:\cB\to \cB$ on a Banach space $\cB$ can be expressed as 
\[
\inf\{\; r(\LL|_{W})\;\mid \mbox{ $W\subset \cB$ is a closed $\LL$-invariant subspace of
finite codimension.}\}\, ,
\]
where $r(\LL|_{W})$ is the spectral radius of the restriction of $\LL$ to $W$. 
Indeed, take any $\tilde{\rho}>r_{ess}(\LL)$, and let $W$ be the image of the spectral projector 
corresponding to the part of spectrum in the disk $\{|z|<\tilde{\rho}\}$, then we see that the infimum 
above is not greater than $\tilde{\rho}$ and hence not greater than $r_{ess}(\LL)$. 
Next let $W$ be an arbitrary closed $\LL$-invariant subspace of 
finite codimension, and let $W'$ be a complementary subspace of $W$ in $\cB$ of finite 
dimension. Let $\pi:\cB\to W$ and $\pi':\cB\to W'$ be the projections corresponding 
to the decomposition $\cB=W\oplus W'$. Then we can decompose 
$\LL$ as $\LL=\LL\circ \pi +\LL\circ \pi'$, where $\LL\circ \pi'$ is of finite rank. 
This implies that the essential spectral radius of 
$\LL$ is bounded by $r(\LL\circ \pi)=r(\LL|_W)$ and hence by  the infimum above. 

The intersection $B_1\cap B_2$ is a Banach space with respect to the norm 
$\|\cdot\|_1+\|\cdot \|_2$. 
From the definition above, we can see that the essential spectral radius of the restriction 
$L:B_1\cap B_2\to B_1\cap B_2$ is bounded  by the maximum of those of $L:B_1\to B_1$ and 
$L:B_2\to B_2$. Thus, to prove the lemma,  we may and do assume 
$B_1\subset B_2$ and $\|\cdot \|_2\le \|\cdot\|_1$. 

Consider $\rho>0$ as in the statement of the lemma. Let $E\subset B_1$ be the finite dimensional subspace that is the 
sum of generalized eigenspaces of $L:B_1\to B_1$ for eigenvalues in 
$\{z\in \complex\mid |z|\ge  \rho\}$. Replacing $B_1$ and $B_2$ by their factor space by 
$E$ respectively, we may and do assume that $E=\{0\}$ or that the spectral radius of $L:B_1\to B_1$ is strictly smaller than $\rho$. 

We can now complete the proof by showing that 
$L:B_2\to B_2$ has no eigenvalues greater than or equal to 
$\rho$ in absolute value. Suppose that it were not true. Then we could take an eigenvector 
 for $L:B_2\to B_2$ corresponding to 
an eigenvalue $\lambda$ so that $|\lambda|$ is equal to the spectral radius of $L:B_2\to B_2$ and 
is not less than $\rho$. 
Since $B_0$ is dense in $B_2$, this would imply that there exists a vector $v\in B_0\subset B_1$ such that 
$\|L^n v\|_1\ge \|L^n v\|_2 \ge |\lambda|^n$ for all $n\ge 0$. 
This contradicts the fact that the spectral radius of $L:B_1\to B_1$ is
strictly smaller than $\rho$. 
\end{proof}

Since the spaces of functions in this paper and \cite{BT}, as well as those in \cite{GoL}, are 
completions of the space of $C^{r-1}$ functions and embedded in the space of 
distributions, the lemma above tells that the part of spectrum of $\LL$ outside of the essential 
spectral radius does not depend on those choices of  Banach spaces.

In view of Lemma~\ref{changespace}, it is natural to
ask whether there exists a Banach space containing
$C^{r-1}(V)$ on which $\LL$ is bounded and has essential
spectral radius strictly smaller than $Q_{r-1}(T,g)$.
We expect that 
there may be such Banach spaces if $d\ge 3$ but not if $d=2$.  
(For hyperbolic endomorphisms, we can find examples of such Banach spaces in \cite[Thereom 3]{AGT}.)

\section{Proof of the inequality (\ref{ineqq})}
\label{ineq}

We show here the claim in Remark~\ref{ineqq} that
$\rho^{p,q}(T,g)\le \inf_{t\in [1, \infty]}R^{p,q,t}(T,g)$. The proof will imply
that the inequality can be strict.
Put
\[
R^{p,q,t}(T,g,m)= \sup_{X} |\det DT^m|^{-1/t}(x) |g^{(m)}(x)| \lambda^{(p,q,m)}(x)
 \]
  for $m\ge 1$ and $t\in [1,\infty]$. Since $\supp(g)$ is contained in $V$, we have
\begin{equation}\label{eqn:RR}
R^{p,q,t}(T,g)=\lim_{m\to \infty}\left(
R^{p,q,t}(T,g,m) \right)^{1/m} \, .
\end{equation}
For each $m\ge 0$, we put $V^m=\cap_{k=0}^{m-1} T^{-k}(V)$, $V^m_+=\{x\in V^m;  |\det DT^m(x)|\ge 1\}$ and 
$V^m_-=V^m-V^m_+$. Then we have
\begin{align*}
\int_X (|g^{(m)}|\cdot \lambda^{(p,q,m)})(x) dx &\le \int_{V^m_-} (|g^{(m)}|\cdot \lambda^{(p,q,m)})(x) dx
\\
&\quad + \int_{T^m(V^m_+)} \!\!\!(|\det DT^m|^{-1}\cdot |g^{(m)}| \lambda^{(p,q,m)})(T^{-m}(y)) dy \, .
\end{align*}
Since  $|\det DT^m|^{-1/t}\ge 1$ on $V^m_-$ and $|\det DT^m|^{-1/t}\ge |\det DT^m|^{-1}$ on $V^m_+$, and since
$V^m_-\subset V$ and $T^m(V^m_+)\subset T(V)$,  we obtain
\[
\int_X (|g^{(m)}|\cdot \lambda^{(p,q,m)})(x) dx\le C \cdot R^{p,q,t}(T,g,m)
\]
for any $m\ge 1$ and $t\in [1,\infty]$. Now the claim follows from (\ref{eqn:RR}) and the definition of $\rho^{p,q}(T,g)$.



\section{The key estimate  (\ref{inte}) for integer $r$}
\label{apd:r}

For the convenience of the reader, we reproduce here the proof of
(\ref{inte}) from \cite{BT}, when $r\ge 2$ is an integer.

By "integration by parts on $w$", we will mean application, 
for $f\in C^2(\real^d)$ and  $g\in C^1_0(\real^d)$,
with $\sum_{j=1}^{d}(\partial_j f)^2\ne 0$ on $\supp(g)$,
of the formula\footnote{To handle noninteger $r>1$, we may use "regularised integration by parts" instead. See \cite{BT}.}
\begin{align*}
 \int_{\real^d} e^{if(w)}g(w) \, dw
&=i\cdot \int_{\real^d} e^{if(w)}\cdot \sum_{k=1}^d \partial_k\left(\frac
{\partial_k f(w)\cdot  g(w)}
{\sum_{j=1}^{d}(\partial_j f(w))^2}\right) \, dw\,  .
\end{align*}
Rewrite the operator $S^{\ell,\tau}_{n,\sigma}$ as 
$
S^{\ell,\tau}_{n,\sigma}(u)(x)=\int_{\real^d} V_{n,\sigma}^{\ell,\tau}(x,y)|\det D\TT| u(\TT(y))\, dy
$,
where
\begin{equation}\label{Vkernel}
V_{n,\sigma}^{\ell,\tau}(x,y)=(2\pi)^{-2d}\int e^{i(x-w)\xi+i(\TT(w)-\TT(y))\eta} 
G(w)\psi_{\Theta',n,\sigma}(\xi)
\tilde{\psi}_{\Theta,\ell,\tau}(\eta)\, dw \, d\xi d\eta \, .
\end{equation}
The required estimate thus  follows if we show, for some  $C(\TT,G)>0$, that 
\begin{equation}\label{eqn:Kernelest}
|V_{n,\sigma}^{\ell,\tau}(x,y)|\le C(\TT,G) 2^{-(r-1)\max\{n,\ell\}}\cdot  
2^{d\min\{n,\ell\}}b(2^{\min\{n,\ell\}}(x-y))
\end{equation}
for all $(\ell,\tau), (n,\sigma) \in \Gamma$ with $(\ell,\tau)\not\hookrightarrow (n,\sigma)$.
Recall the constant $C(\TT)$ in (\ref{lowerbd}).
W may and do assume  $\max\{n,\ell\}> C(\TT)$ in proving 
(\ref{eqn:Kernelest}).
Integrating (\ref{Vkernel}) by parts $(r-1)$ times on $w$, we obtain 
\[
V_{n,\sigma}^{\ell,\tau}(x,y)=(2\pi)^{-2d}
\int e^{i(x-w)\xi+i(\TT(w)-\TT(y))\eta} 
F(\xi,\eta,w)\psi_{\Theta',n,\sigma}(\xi)
\tilde{\psi}_{\Theta,\ell,\tau}(\eta)\, dw\, d\xi\,  d\eta 
\]
where $F(\xi,\eta,w)$ is $C^\infty$ in the variables  $\xi$ and $\eta$,  continuous in $w$ and supported on
$\real^d\times \real^d\times \supp(G)$.
Using (\ref{lowerbd}), we can see that, 
if $\psi_{\Theta',n,\sigma}(\xi)
\tilde{\psi}_{\Theta,\ell,\tau}(\eta)\ne 0$, then, for each multi-indices $\alpha$, $\beta$ and for $\xi,\eta\in \real^d$, 
\begin{equation}\label{exest}
\|\partial_\xi^\alpha\partial_\eta^\beta
F(\xi,\eta,\cdot)\|_{C^{0}}\le C_{\alpha,\beta}(\TT,G)
2^{-n|\alpha|-\ell|\beta|-(r-1)\max\{n,\ell\}}\, .
\end{equation}
Put 
$H_{n,\ell}(\xi,\eta,w)=F(\xi,\eta,w)\psi_{\Theta',n,\sigma}(\xi)
\tilde{\psi}_{\Theta,\ell,\tau}(\eta)
$,
and consider the scaling
\[
\widetilde H_{n,\ell}(\xi,\eta,w)=H_{n,\ell}(2^{n}\xi, 2^{\ell}\eta,w)\, .
\]
The estimate (\ref{exest}) implies that for all $\alpha$ and $\beta$
\[
\|\partial_\xi^\alpha\partial_\eta^\beta
\widetilde H_{n,\ell}(\xi,\eta,\cdot)\|_{C^{0}}\le C_{\alpha,\beta}(\TT,G)
 2^{-(r-1)\max\{n,\ell\}}\, .
\]
Denote by $\Fourier^{-1}_{\xi\eta}$ the inverse Fourier transform with respect to the variables $\xi$ and $\eta$. Then the estimate above on $\widetilde H_{n,\ell}$ implies 
\begin{align*}
|(\Fourier^{-1}_{\xi\eta}H_{n,\ell})(x,y,w)|
&=|(\Fourier^{-1}_{\xi\eta} \widetilde{H}_{n,\ell})(2^{n}x,2^{\ell}y,w)|\\
&\le C_{\alpha,\beta}(\TT,G)
\cdot 2^{-(r-1)\max\{n,\ell\}} \cdot b_{n}(x) \cdot b_{\ell}(y)
\end{align*}
where $b_{m}$ is the function defined in (\ref{eqn:bm}). 
Therefore we obtain  
\begin{align*}
|V^{n,\sigma}_{\ell,\tau}(x,y)|&=\left|\int_{\supp(G)} (\Fourier^{-1}_{\xi\eta}H_{n,\ell})(x-w,\TT(w)-\TT(y),w) dw \right|\\
&\le C(\TT,G) \cdot 2^{-(r-1)\max\{n,\ell\}} \cdot 
\int b_{n}(x-w) \cdot b_{\ell}(w-y) dw\,, 
\end{align*}
where we used the fact that $\TT$ is bilipschitz to replace 
$\TT(w)-\TT(y)$ by $w-y$. Now (\ref{eqn:Kernelest}) follows from (\ref{eqn:bconv}). 

\bibliographystyle{amsplain}

\end{document}